\documentclass{article}

\usepackage{arxiv}
\usepackage{soul}
\usepackage{comment}
\usepackage{xcolor}
\usepackage[utf8]{inputenc} 
\usepackage[T1]{fontenc}    
\usepackage{hyperref}       
\usepackage{url}            
\usepackage{booktabs}       
\usepackage{amsfonts}       
\usepackage{nicefrac}       
\usepackage{microtype}      
\usepackage{lipsum}

\usepackage{morefloats}
\usepackage{color}
\usepackage{multirow}
\usepackage{latexsym}
\usepackage{amsmath,amssymb,amsthm,graphicx}
\usepackage{dsfont}
\usepackage{soul}  
\usepackage{enumitem}
\usepackage{hyperref}
\usepackage{arydshln}
\usepackage{lscape}
\usepackage[normalem]{ulem}
\newtheoremstyle{thm}
{9pt}
{9pt}
{\itshape}
{}
{\bfseries}
{.}
{ }
{}
\theoremstyle{thm}
\DeclareMathAlphabet{\mathscr}{OT1}{pzc}{m}{it}
\newtheorem{theorem}{Theorem}[section]

\newtheoremstyle{def}
{9pt}
{9pt}
{}
{}
{\bfseries}
{.}
{ }
{}
\theoremstyle{def}

\newtheorem{example}[theorem]{Example}

\newcommand{\R}{\mathbb{R}} 
\newcommand{\N}{\mathbb{N}} 
\newcommand{\E}{\mathbb{E}} 

\renewcommand{\footnoterule}{%
	\kern -3.5pt
	\hrule width \textwidth height 1pt
	\kern 3.5pt
}

\newcommand{\be}[1]{\textcolor{red}{[BE: #1]}}

\makeatletter
\def\blfootnote{\xdef\@thefnmark{}\@footnotetext}
\makeatother

\title{Eigenvalues approximation of integral covariance operators with applications to weighted $L^2$ statistics}

\author{Bruno Ebner\\
 Institute of Stochastics, \\
Karlsruhe Institute of Technology (KIT), \\
Englerstr. 2, D-76133 Karlsruhe, Germany. \\
\href{mailto:Bruno.Ebner@kit.edu}{Bruno.Ebner@kit.edu}
\And
M. Dolores Jim\'{e}nez-Gamero\\
Dpto. Estadística e Investigación Operativa,\\
Universidad de Sevilla,\\
Avda Reina Mercedes sn,\\
41012-Sevilla, Spain\\
\href{mailto:dolores@us.es}{dolores@us.es}
\And
Bojana Milo\v{s}evi\'{c}\\
University of Belgrade\\
Faculty of Mathematics\\
\href{mailto:bojana@matf.bg.ac.rs}{bojana@matf.bg.ac.rs}
}

\begin{document}

\date{\today}

\maketitle

\begin{abstract}
	Finding the eigenvalues connected to the covariance operator of a centred Hilbert-space valued Gaussian process is genuinely considered a hard problem in several mathematical disciplines. In statistics this problem arises for instance in the asymptotic null distribution of goodness-of-fit test statistics of weighted $L^2$-type. For this problem we present the Rayleigh-Ritz method to approximate the eigenvalues. The usefulness of these approximations is shown by high lightening implications such as critical value approximation and  theoretical comparison of test statistics by means of Bahadur efficiencies.
\end{abstract}

\section{Introduction}\label{sec:Intro}
In recent decades, numerous goodness-of-fit tests have been developed based on statistics of weighted $L^2$-type. These methods rely on empirical characteristic functionals of the distribution, such as the empirical distribution function, the empirical characteristic function, the empirical Laplace transform, the empirical moment generating function, the empirical Hankel transform, or empirical Stein-type characterizations, to assess deviations from their theoretical counterparts or zero. For a comprehensive list of weighted $L^2$ statistics up to 2015, we refer to \cite{BEH:2017}. Our focus, however, is on more recent advancements, such as those discussed in \cite{AES:2023,AP:2017,BBQ:2018,BE:2019b,HR:2020b,HJ:2019,HJM:2019,HM:2019,HNNS:2019}.

Theoretical results on weighted $L^2$-statistics often involve a nondegenerate limit distribution of the statistic (say) $T_n$ under the null hypothesis, where the authors typically demonstrate that $T_n$ converges weakly to $\|Z\|^2$, where $Z$ is a Gaussian process taking values in a suitable function space and $\|\cdot\|$ is a corresponding norm. Most authors then refer to the renowned Karhunen-Loeve transform, showing the equality in distribution of $\|Z\|^2$ and $\sum_{j=1}^\infty \lambda_j N_j^2$, where $N_1, N_2, \ldots$ are independent identically distributed (iid) standard normal random variables, and $\{\lambda_j\}_{j\in\N}$ is the sequence of positive eigenvalues of the covariance operator of $Z$. Since the involved covariance kernels of the Gaussian process are often explicitly known but exhibit a certain complexity, further calculations are typically halted, and the corresponding critical values for performing the test are determined via Monte Carlo simulation methods. In some very special cases, the eigenvalues of the covariance operators can be calculated analytically, see \cite{BGT:2018,BH:2008,BL:2005,E:2023,HN:2000,MT:2005}, all of which are related to the classical eigenvalue problem of the Brownian bridge covariance kernel. For a rigorous derivation of the solution, see \cite{DM:2003,SW:1986}. Some articles address the approximation of the eigenvalues using the quadrature method, see \cite{BL:2005}, by finding the roots of the connected Fredholm determinant numerically, see \cite{EH:2023,S:1976,S:1977}, by employing a stochastic Monte Carlo type approach, see \cite{EH:2023b,FanEtAl2017}, or by approximating the kernel \cite{bozin2020new, JAJB:2015, NJ:2016}.

The primary objective of this paper is to utilize the Rayleigh-Ritz method (previously overlooked in statistical literature) to approximate the largest $m$ eigenvalues $\lambda_1, \dots, \lambda_m$ of the covariance kernel, as referenced in \cite{K:1971}, Section 7.9. To this end, the paper is organized as follows. Section \ref{sec:EVP} presents the eigenvalue problem, while Section \ref{RR-method} details the Rayleigh-Ritz method. Section \ref{sec:AppL2} applies the method to both established and novel eigenvalue problems in the literature for various supports. Finally, Section \ref{sec:Con} discusses the implications of the knowledge of the eigenvalues, including a method for approximating the quantiles of the asymptotic distribution and providing approximate Bahadur efficiencies for several examples.

\section{The eigenvalue problem}\label{sec:EVP}
Using the notation from \cite{BEH:2017}, let $\mathcal{B}^d$ denote the Borel $\sigma$-algebra of subsets of $\R^d$, with $M \neq \emptyset$ being an element of $\mathcal{B}^d$ and $\mu$ a finite measure on $\mathcal{B}^d_M = \mathcal{B}^d \cap M$. Define $L^2 = L^2(M, \mathcal{B}_M^d, \mu)$ as the (separable) Hilbert space consisting of (equivalence classes of) square-integrable measurable functions on $M$, equipped with the inner product $\langle g, h \rangle = \int_M gh \, \mathrm{d} \mu$. Moreover, let $\|h\|^2_{L^2} = \int_M h^2 \, \mathrm{d} \mu$ and let $Z = (Z(t), t \in M)$ be a centered Gaussian process that can be considered a random element of $L^2$ with $\E \|Z\|^2_{L^2} < \infty$. The distribution of $Z$ is uniquely characterized by its covariance kernel $K(s,t)$, $s,t \in M$. Consequently, as well-known (see \cite{DW:2014}, Chapter 3), the distribution of $\|Z\|^2_{L^2}$ corresponds to $\sum_{j=1}^\infty \lambda_j N_j^2$, where $N_j$ are iid standard normal random variables and $\lambda_j$ are the positive eigenvalues, at most countable, associated with eigenfunctions $f_j$ of the (linear second-order homogeneous Fredholm) integral equation
\begin{equation}\label{int:eq}
\lambda f(s) = \mathcal{K} f(s), \quad s \in M,
\end{equation}
where $\mathcal{K}: L^2 \mapsto L^2$ is defined as
\[
\mathcal{K} f(s) = \int_M K(s,t) f(t) \mu(\mathrm{d} t),
\]
corresponding to the covariance kernel $K$ of $Z$ (for properties of covariance kernels, see \cite{SW:1986}, p.207). The task of finding eigenvalues and eigenfunctions is often referred to as the kernel eigenproblem; see \cite{WS:2003}. Solving \eqref{int:eq} analytically is generally considered challenging. For a list of solutions for specific choices of $\mathcal{K}$, $M$, and $\mu$, see Appendix A in \cite{FM:2016}. From now on, we assume that $\mu$ is defined by a non-negative weight function $w$ on $M$, such that $\mu(\mathrm{d} t) = w(t) \mathrm{d} t$ and $\int_M t w(t) \mathrm{d} t < \infty$. We will also consider an example where $M = \mathbb{N}_0 = \{0, 1, 2, \ldots\}$ and $\mu$ has a density (in the Radon-Nikodym sense) with respect to the counting measure.

\section{The Rayleigh-Ritz method} \label{RR-method}
As mentioned in the Introduction, let $K(s,r)$ be a covariance kernel. Given that it is symmetric and nonnegative definite, its eigenvalues are real and nonnegative. Throughout this paper, we assume implicitly that $\int_M K(t,t) \mu(\mathrm{d} t) < \infty$, which in turn implies that $\int_M \int_M K(s,t)^2 \mu(\mathrm{d} s) \mu(\mathrm{d} t) < \infty$.

For the sake of completeness, this section reformulates the Rayleigh-Ritz method. Let $\lambda_1$ denote the largest eigenvalue of the covariance kernel $K$ and let $f_1$ be the associated normalized eigenfunction, such that $\|f_1\|_{L^2}^2=1$. Let $\{\psi_j\}_{j \geq 1}$ represent an orthonormal basis of $L^2$. The core idea is to approximate $f_1$ by $f_{1n} = \sum_{j=1}^n \alpha_j \psi_j$, with $\alpha_1, \ldots, \alpha_n \in \R$, for some $n \in \N$. In this approximation, $\alpha_j = \langle f_1, \psi_j \rangle$, for $1 \leq j \leq n$, which are unknown quantities (as $f_1$ is unknown). To estimate $\lambda_1=\langle \mathcal{K} f_{1}, f_{1} \rangle$,  we maximize the function
\begin{equation*}
\langle \mathcal{K} f_{1n}, f_{1n} \rangle = \langle K \sum_{j=1}^n \alpha_j \psi_j, \sum_{j=1}^n \alpha_j \psi_j \rangle = \sum_{j,k=1}^n K_{jk} \alpha_j \alpha_k,
\end{equation*}
subject to
\begin{equation*}
\left\| \sum_{j=1}^n \alpha_j \psi_j \right\|_{L^2}^2 = \sum_{j=1}^n \alpha_j^2 = 1,
\end{equation*}
where $K_{jk} = \langle K \psi_j, \psi_k \rangle = \langle \psi_j, K \psi_k \rangle$. Using the method of Lagrangian multipliers, we set
\begin{equation*}
\mathbf{L}(x) = \sum_{j,k=1}^n K_{jk} \alpha_j \alpha_k - x \delta_{jk} \alpha_j \alpha_k, \quad x \in \R,
\end{equation*}
where $\delta_{jk}$ is the Kronecker delta. The extremal values of $\alpha_1, \ldots, \alpha_k$ are determined from the equations $\partial \mathbf{L}(x) / \partial \alpha_j = 0$ for $j = 1, \ldots, n$, leading to
\begin{equation*}
\sum_{k=1}^n K_{jk} \alpha_k - x \alpha_j = 0, \quad j = 1, \ldots, n.
\end{equation*}
This system of equations has a nontrivial solution if and only if the determinant
\begin{equation}\label{eq:det}
\det \begin{pmatrix}
K_{11} - x & K_{12} & \cdots & K_{1n} \\
\vdots & \vdots & \ddots & \vdots \\
K_{n1} & K_{n2} & \cdots & K_{nn} - x
\end{pmatrix} = 0.
\end{equation}
Thus, the method calculates the eigenvalues of the $n \times n$ matrix $M_n = (K_{jk})_{1 \leq j,k \leq n}$. For $1 \leq m \leq n$, the $m$ largest eigenvalues of $M_n$, denoted by $\widehat{\lambda}_1, \ldots, \widehat{\lambda}_m$, estimate the $m$ largest eigenvalues of the covariance kernel $K$, denoted by $\lambda_1, \ldots, \lambda_m$. Larger values of $n$ will improve the approximation. In fact, Theorem \ref{thm} below shows that $|\widehat{\lambda}_i - \lambda_i| \to 0$ as $n \to \infty$. Before stating and proving this result, we first see that the Rayleigh-Ritz method approximates not only $f_1$ by $f_{1n}$, but also the kernel $K$.

If $\{\psi_j(s), \, s \in M \}_{j \geq 1}$ is an orthonormal basis of $L^2$, then $\{\psi_j(s) \psi_k(t), \, s,t \in M \}_{j,k \geq 1}$ forms an orthonormal set in the (separable) Hilbert space ${\cal H} = \{h : M \times M \to \mathbb{R} \, | \, \|h\|^2_{\cal H} = \int_M \int_M h(s,t)^2 \mu(\mathrm{d}s) \mu(\mathrm{d}t) < \infty\}$, with the inner product
$\langle h_1, h_2 \rangle_{\cal H} = \int_M \int_M h_1(s,t) h_2(s,t) \mu(\mathrm{d}s) \mu(\mathrm{d}t)$, for $h_1, h_2 \in {\cal H}$. Thus, one can approximate the covariance kernel $K(s,t)$ by
\begin{equation} \label{Kn}
K_n(s,t) = \sum_{j,k=1}^n K_{jk} \psi_j(s) \psi_k(t).
\end{equation}
It is easy to see that $K_n$ is symmetric and nonnegative definite, hence its eigenvalues are real and nonnegative. Moreover, it has at most $n$ positive eigenvalues. Let ${\cal K}_n$ denote the operator analogous to ${\cal K}$, defined using $K_n$ instead of $K$. It is straightforward to verify that
\begin{equation*}
\langle {\cal K} f_{1n}, f_{1n} \rangle = \langle {\cal K}_n f_{1}, f_{1} \rangle = \langle {\cal K}_n f_{1n}, f_{1n} \rangle.
\end{equation*}
Therefore, the Rayleigh-Ritz method essentially approximates both $f_1$ by $f_{1n}$ and $K$ by $K_n$. This observation is crucial for proving the following result.

\begin{theorem} \label{thm}
Let $\{\lambda_i\}_{i \geq 1}$ denote the eigenvalues of $K$, arranged in decreasing order and repeated according to their multiplicity. Let $\{\widehat{\lambda}_i\}_{i=1}^n$ represent the solutions (in $x$) of equation \eqref{eq:det}, also arranged in decreasing order and repeated according to their multiplicity. Then, $|\widehat{\lambda}_i - \lambda_i| \to 0$ as $n \to \infty$, for each $i$.
\end{theorem}

{\sc Proof}
Section 7.7 of \cite{K:1971}, shows that
\begin{equation} \label{KKn}
\|K-K_n\|_{\cal H} \to 0,\quad \mbox{as $n\to \infty$},
\end{equation}
where $K_n$ is as defined in \eqref{Kn}.
The result follows from \eqref{KKn} and  \cite{DSpartII}, p. 1090-1091. $\Box$

\medskip

It is evident that the choice of the orthonormal set $\{\psi_j\}_{j \geq 1}$ depends on the support $M$ and the weight function $w(t)$. We provide examples for the most common supports and weight functions, utilizing orthogonal polynomials. For a list of classical polynomials, see Table 1.1 on p. 29 of \cite{G:2004}.


\section{Applications} 
\label{sec:AppL2}
In this section, we apply the Rayleigh-Ritz method to both solved and open eigenvalue problems found in the literature. The former is used to demonstrate the quality of the approximation, while the latter showcases the method's flexibility. We provide appropriate orthonormal bases associated with various supports and weight functions of the Hilbert spaces $L^2$. Hereafter, we denote $u \wedge v = \min(u,v)$ and $u \vee v = \max(u,v)$, for $u,v \in \mathbb{R}$.

\subsection{Support $M=[0,1]$}
If the domain of integration is denoted by $M=[0,1]$, we concentrate on the weight function $w(t)=1$. A suitable set of orthonormal polynomials can then be expressed as
\begin{equation*}
    \phi_k(x)=\sqrt{2k+1}P_k(2x-1),\quad x\in[0,1],
\end{equation*}
where
$ 
    P_k(x)=\sum_{j=0}^{\lfloor k/2\rfloor}(-1)^j\frac{(2k-2j)!}{(k-j)!(k-2j)!j!2^k}x^{k-2j}
$ 
represents the Legendre polynomial of degree $k$.

\begin{example} (Non-parametric goodness-of-fit tests)
\begin{itemize}
\item The traditional Cram\'{e}r-von Mises test statistic possesses the covariance kernel
\begin{equation} \label{CvM}
    K(s,t)=s\wedge t - st, \quad s,t\in[0,1],
\end{equation}
refer to \cite{SW:1986}, p. 14 (see also \cite{M:2020}). It is a well-established fact that the sequence of eigenvalues is given by $\lambda_j=1/j^2\pi^2$, $j=1,2,\ldots,$
Table \ref{tab:CvM} presents the estimators of the five largest eigenvalues of kernel \eqref{CvM}, for $3 \leq n \leq 15$, along with the actual values.  Throughout this paper, e-$j$ stands for $10^{-j}$. Observing this table, we notice that, with the precision in it, for $n \geq 4$ the estimator of the largest eigenvalue matches the actual value; $n \geq 9$ for the second; $n \geq 10$ for the third, and so forth.

\begin{table}[t]
    \centering
    \begin{tabular}{rrrrrr}
    $n$ &  \multicolumn{1}{c}{$\widehat{\lambda}_1$} &  \multicolumn{1}{c}{$\widehat{\lambda}_2$} &  \multicolumn{1}{c}{$\widehat{\lambda}_3$} &
           \multicolumn{1}{c}{$\widehat{\lambda}_4$} &  \multicolumn{1}{c}{$\widehat{\lambda}_5$}\\ \hline
    3 &  1.012648e-1 &  2.514783e-2 &  0.5878034e-2 &  0.000000e-3 &  0.000000e-3\\
    4 &  1.013212e-1 &  2.514783e-2 &  1.0961853e-2 &  2.629952e-3 &  0.000000e-3\\
    5 &  1.013212e-1 &  2.532945e-2 &  1.0961853e-2 &  5.955316e-3 &  1.353354e-3\\
    5 &  1.013212e-1 &  2.532945e-2 &  1.0961853e-2 &  5.955316e-3 &  1.353354e-3\\
    6 &  1.013212e-1 &  2.532945e-2 &  1.1253177e-2 &  5.955316e-3 &  3.626214e-3\\
    7 &  1.013212e-1 &  2.533029e-2 &  1.1253177e-2 &  6.318910e-3 &  3.626214e-3\\
    8 &  1.013212e-1 &  2.533029e-2 &  1.1257885e-2 &  6.318910e-3 &  4.025014e-3\\
    9 &  1.013212e-1 &  2.533030e-2 &  1.1257885e-2 &  6.332401e-3 &  4.025014e-3\\
   10 &  1.013212e-1 &  2.533030e-2 &  1.1257909e-2 &  6.332401e-3 &  4.052147e-3\\
   11 &  1.013212e-1 &  2.533030e-2 &  1.1257909e-2 &  6.332573e-3 &  4.052147e-3\\
   12 &  1.013212e-1 &  2.533030e-2 &  1.1257909e-2 &  6.332573e-3 &  4.052840e-3\\
   13 &  1.013212e-1 &  2.533030e-2 &  1.1257909e-2 &  6.332574e-3 &  4.052840e-3\\
   14 &  1.013212e-1 &  2.533030e-2 &  1.1257909e-2 &  6.332574e-3 &  4.052847e-3\\
   15 &  1.013212e-1 &  2.533030e-2 &  1.1257909e-2 &  6.332574e-3 &  4.052847e-3\\
 true &  1.013212e-1 &  2.533030e-2 &  1.1257909e-2 &  6.332574e-3 &  4.052847e-3
 \end{tabular}
 \medskip
    \caption{Estimators of the five largest eigenvalues of kernel \eqref{CvM} obtained for $3 \leq n \leq 15$  with the Rayleigh-Ritz method.} 
    \label{tab:CvM}
\end{table}

\item In \cite{HN:2000}, equation (3.2), the covariance kernel is characterized by
\begin{equation} \label{CvM2}
    \overline{K}_0(s,t)=\frac{st(s\wedge t)}{2}-\frac{(s\wedge t)^3}{6}-\frac{s^2t^2}{4},\quad s,t\in[0,1].
\end{equation}
Table \ref{tab:CvM2} presents the estimators of the five largest eigenvalues of kernel \eqref{CvM2}, for $3 \leq n \leq 15$, along with the actual values, which are taken from Table 3.1 of \cite{HN:2000}. Observing this table, we notice that, with the precision in it, for $n \geq 5$ the estimator of the largest eigenvalue matches the actual value with exception of the third eigenvalue marked by $^\ast$. In this case we conjecture that the authors of \cite{HN:2000} made a numerical or typographical error in the fourth decimal place.

\begin{table}[t]
    \centering
    \begin{tabular}{rrrrrr}
    $n$ &  \multicolumn{1}{c}{$\widehat{\lambda}_1$} &  \multicolumn{1}{c}{$\widehat{\lambda}_2$} &  \multicolumn{1}{c}{$\widehat{\lambda}_3$} &
           \multicolumn{1}{c}{$\widehat{\lambda}_4$} &  \multicolumn{1}{c}{$\widehat{\lambda}_5$}\\ \hline
 3 & 3.196304e-2 & 1.058647e-3 & 1.564775e-4 & 0.000000e-5 & 0.000000e-5\\
 4 & 3.196394e-2 & 1.093812e-3 & 1.633718e-4 & 3.999208e-5 & 0.000000e-5\\
 5 & 3.196395e-2 & 1.094444e-3 & 1.787795e-4 & 4.309337e-5 & 1.330054e-5\\
 6 & 3.196395e-2 & 1.094568e-3 & 1.792433e-4 & 5.121087e-5 & 1.484606e-5\\
 7 & 3.196395e-2 & 1.094568e-3 & 1.794972e-4 & 5.155691e-5 & 1.950214e-5\\
 8 & 3.196395e-2 & 1.094569e-3 & 1.795006e-4 & 5.189992e-5 & 1.976122e-5\\
 9 & 3.196395e-2 & 1.094569e-3 & 1.795022e-4 & 5.190770e-5 & 2.014210e-5\\
10 & 3.196395e-2 & 1.094569e-3 & 1.795022e-4 & 5.191281e-5 & 2.015456e-5\\
11 & 3.196395e-2 & 1.094569e-3 & 1.795022e-4 & 5.191289e-5 & 2.016592e-5\\
12 & 3.196395e-2 & 1.094569e-3 & 1.795022e-4 & 5.191292e-5 & 2.016616e-5\\
13 & 3.196395e-2 & 1.094569e-3 & 1.795022e-4 & 5.191292e-5 & 2.016629e-5\\
14 & 3.196395e-2 & 1.094569e-3 & 1.795022e-4 & 5.191292e-5 & 2.016629e-5\\
15 & 3.196395e-2 & 1.094569e-3 & 1.795022e-4 & 5.191292e-5 & 2.016629e-5\\
true & 3.196395e-2 & 1.094569e-3 & 1.795105e-4$^\ast$ & 5.191292e-5 & 2.016629e-5
 \end{tabular} \medskip
    \caption{Estimators of the five largest eigenvalues of kernel \eqref{CvM2} obtained for $3 \leq n \leq 15$  with the Rayleigh-Ritz method.}
    \label{tab:CvM2}
\end{table}

\item In \cite{ebner2021new}, Theorem 2.1, the covariance kernel is defined as
\begin{equation} \label{CvM3}
K_{Z}(s,t)=\frac{1-(2(s \vee t)-1)^3}{6}-st(1-s)(1-t),\quad s,t\in(0,1).
\end{equation}
Table \ref{tab:CvM3} presents the estimators of the five largest eigenvalues of kernel \eqref{CvM3}, for $3 \leq n \leq 15$. Similar to the previous two cases, we notice a quick convergence of the first few eigenvalues.

\begin{table}[t]
    \centering
    \begin{tabular}{rrrrrr}
    $n$ &  \multicolumn{1}{c}{$\widehat{\lambda}_1$} &  \multicolumn{1}{c}{$\widehat{\lambda}_2$} &  \multicolumn{1}{c}{$\widehat{\lambda}_3$} &
           \multicolumn{1}{c}{$\widehat{\lambda}_4$} &  \multicolumn{1}{c}{$\widehat{\lambda}_5$}\\ \hline
 3 & 1.155505e-1 & 6.173688e-3 & 2.232223e-3 & 0.000000e-3 & 0.0000000e-3\\
 4 & 1.157325e-1 & 6.553968e-3 & 2.287779e-3 & 1.615387e-3 & 0.0000000e-3\\
 5 & 1.157342e-1 & 6.936156e-3 & 2.451756e-3 & 1.636573e-3 & 0.8705469e-3\\
 6 & 1.157342e-1 & 6.941384e-3 & 2.699774e-3 & 1.851195e-3 & 0.8713323e-3\\
 7 & 1.157342e-1 & 6.941476e-3 & 2.721395e-3 & 1.946449e-3 & 1.0884916e-3\\
 8 & 1.157342e-1 & 6.942877e-3 & 2.723814e-3 & 1.971728e-3 & 1.1784099e-3\\
 9 & 1.157342e-1 & 6.943677e-3 & 2.725788e-3 & 1.973082e-3 & 1.1888306e-3\\
10 & 1.157342e-1 & 6.943684e-3 & 2.728232e-3 & 1.976252e-3 & 1.1888330e-3\\
11 & 1.157342e-1 & 6.943684e-3 & 2.728373e-3 & 1.977335e-3 & 1.1955965e-3\\
12 & 1.157342e-1 & 6.943684e-3 & 2.728383e-3 & 1.977537e-3 & 1.1976705e-3\\
13 & 1.157342e-1 & 6.943685e-3 & 2.728387e-3 & 1.977544e-3 & 1.1978559e-3\\
14 & 1.157342e-1 & 6.943685e-3 & 2.728391e-3 & 1.977555e-3 & 1.1978559e-3\\
15 & 1.157342e-1 & 6.943685e-3 & 2.728392e-3 & 1.977558e-3 & 1.1979070e-3
 \end{tabular} \medskip
    \caption{Estimators of the five largest eigenvalues of kernel \eqref{CvM3} obtained for $3 \leq n \leq 15$  with the Rayleigh-Ritz method.}
    \label{tab:CvM3}
\end{table}


\end{itemize}
\end{example}

\subsection{Support $M=[0,\infty)$}
If the domain of integration is denoted by $M=[0,\infty)$, we concentrate on the weight function $w_\gamma(t)=e^{-\gamma t}$ for a positive parameter $\gamma$. A suitable set of orthonormal polynomials can then be expressed as
\begin{equation*}
    \phi_k(x;\gamma)=\sqrt{\gamma}L_k(\gamma x),\quad x\in[0,\infty),
\end{equation*}
where
$ 
    L_k(x)=\sum_{j=0}^{k}\binom{k}{j}\frac{(-1)^j}{j!}x^j
$ 
represents the Laguerre polynomial of degree $k$.

\begin{example}(Testing Exponentiality) \label{exm: exp}
\begin{itemize}
\item In \cite{BH:2008}, the covariance kernel is specified in (7) as
\begin{equation}\label{kernel rho}
    \rho(s,t)=\min(1-e^{-s},1-e^{-t})-(1-e^{-s})(1-e^{-t}), \quad s,t\ge0.
\end{equation}
In this scenario, the formulas for the eigenvalues are explicitly provided, see Theorem 2 in \cite{BH:2008}.
It's worth noting that in \cite{E:2023}, the author demonstrates that these eigenvalues also correspond to the covariance kernel
\begin{equation}\label{eq:K0}
    K_0(s,t)=e^{-s \vee t}-e^{-(s+t)},\quad s,t\ge 0,
\end{equation}
and the first twenty eigenvalues are tabulated in Table 2 in \cite{E:2023}. Table \ref{Laguerre1} presents the estimators of the two largest eigenvalues of kernel \eqref{eq:K0}, for $n=10(5)30$, and $\gamma=1,2$, along with the actual values, which are taken from Table 2 in \cite{E:2023}.

\begin{table}[t]
    \centering
    \begin{tabular}{rrrrrr}
   & \multicolumn{2}{c}{$\gamma=1$} &  & \multicolumn{2}{c}{$\gamma=2$}\\
   \cline{2-3} \cline{5-6}
$n$ &  \multicolumn{1}{c}{$\widehat{\lambda}_1$} &  \multicolumn{1}{c}{$\widehat{\lambda}_2$} &  &
      \multicolumn{1}{c}{$\widehat{\lambda}_1$}  &  \multicolumn{1}{c}{$\widehat{\lambda}_2$}\\ \hline
10   & 1.012749e-1 & 2.346677e-2 & & 5.274782e-2 & 1.210710e-2\\
15   & 1.013091e-1 & 2.517918e-2 & & 5.275176e-2 & 1.218831e-2\\
20   & 1.013168e-1 & 2.527533e-2 & & 5.275292e-2 & 1.220204e-2\\
25   & 1.013204e-1 & 2.529163e-2 & & 5.275300e-2 & 1.221051e-2\\
30   & 1.013211e-1 & 2.531704e-2 & & 5.275301e-2 & 1.221147e-2\\
true & 1.013212e-1 & 2.533030e-2 & & 5.275301e-2 & 1.221201e-2
 \end{tabular} \medskip
    \caption{Estimators of the two largest eigenvalues of kernel \eqref{eq:K0} obtained for $n=10(5)30$ and $\gamma=1,2$  with the Rayleigh-Ritz method.  }
    \label{Laguerre1}
\end{table}

\item In \cite{K:2001}, Theorem 2.1, the covariance kernel is defined as
\begin{equation} \label{K2.Laguerre}
    K(s,t)=(|s-t|+2)e^{-s\vee t}-(s+t+st+2)e^{-(s+t)}, \quad s,t\ge0.
\end{equation}
Table \ref{Laguerre2} presents the estimators of the two largest eigenvalues of kernel \eqref{K2.Laguerre}, for $n=10(5)30$, and $\gamma=0.5,1,1.5$.

\end{itemize}
\end{example}
\begin{table}[t]
    \centering
    \begin{tabular}{rrrrrrrrr}
   & \multicolumn{2}{c}{$\gamma=0.5$} &  & \multicolumn{2}{c}{$\gamma=1$} &  & \multicolumn{2}{c}{$\gamma=1.5$}\\
   \cline{2-3} \cline{5-6} \cline{8-9}
n &  \multicolumn{1}{c}{$\widehat{\lambda}_1$} &  \multicolumn{1}{c}{$\widehat{\lambda}_2$} &  &
      \multicolumn{1}{c}{$\widehat{\lambda}_1$}  &  \multicolumn{1}{c}{$\widehat{\lambda}_2$}
      &  &
      \multicolumn{1}{c}{$\widehat{\lambda}_1$}  &  \multicolumn{1}{c}{$\widehat{\lambda}_2$}\\ \hline
10 & 6.905560e-2 & 0.9839087e-2 & & 3.091694e-2 & 4.061322e-3 & & 1.581485e-2 & 1.958083e-3 \\
15 & 6.949381e-2 & 1.0149090e-2 & & 3.094987e-2 & 4.154889e-3 & & 1.581630e-2 & 1.980193e-3 \\
20 & 6.963910e-2 & 1.0301102e-2 & & 3.095009e-2 & 4.168487e-3 & & 1.581631e-2 & 1.980447e-3 \\
25 & 6.965260e-2 & 1.0420650e-2 & & 3.095011e-2 & 4.168621e-3 & & 1.581631e-2 & 1.980510e-3 \\
30 & 6.965294e-2 & 1.0445512e-2 & & 3.095012e-2 & 4.168761e-3 & & 1.581631e-2 & 1.980516e-3
 \end{tabular} \medskip
    \caption{Estimators of the two largest eigenvalues of kernel \eqref{K2.Laguerre} obtained for $n=10(5)30$ and $\gamma=0.5,1,1.5$  with the Rayleigh-Ritz method.}
    \label{Laguerre2}
\end{table}

\subsection{Support $M=\R$}
If the domain of integration is denoted by $M=\R$, we concentrate on the weight function $w_\gamma(t)=e^{-\gamma t^2}$ for a positive parameter $\gamma>0$. A suitable set of orthonormal polynomials can then be expressed as
\begin{equation}\label{eq:Hermite}
\phi_k(x;\gamma)=\left(2^k\,k!\sqrt{\pi/\gamma}\right)^{-1/2}H_k\left(\sqrt{\gamma}x\right), \quad x\in\R,
\end{equation}
where $H_k(x)=(-1)^k\exp(x^2)\frac{\mbox{d}^k}{\mbox{d}x^k}\exp(-x^2)$ represents the Hermite polynomial of degree $k$.

\begin{example}(Normality Test) \label{exm: normality}
\begin{itemize}
    \item In \cite{E:2020}, Theorem 2.2, the covariance kernel is defined as
    \begin{equation} \label{K_Z}
        K_Z(s,t)=(st+1)\exp\left(-\frac{(s-t)^2}2\right)-(2st+1)\exp\left(-\frac{s^2+t^2}2\right),\quad s,t\in\R.
    \end{equation}
    The exact values of the four cumulants of the distribution of $T(\gamma)=\|Z\|^2_{L^2}$ with $w=\varphi_\gamma$ (denoted as $\kappa_1(\gamma)=\mathbb{E}(T(\gamma))$, $\kappa_2(\gamma)=\mathbb{E}(T(\gamma)-\kappa_1(\gamma))^2$,  $\kappa_3(\gamma)=\mathbb{E}(T(\gamma)-\kappa_1(\gamma))^3$ and $\kappa_4(\gamma)=\mathbb{E}(T(\gamma)-\kappa_1(\gamma))^4-3\mathbb{E}^2(T(\gamma)-\kappa_1(\gamma))^2$)  have been provided in  \cite{E:2020} for several values of $\gamma$. As seen in Subsection \ref{Pearson.System},
$$
\begin{array}{rclcrcl}
\kappa_1(\gamma) & = & \displaystyle \sum_{j=1}^\infty \lambda_j(\gamma), & &
\kappa_2(\gamma) & = & \displaystyle
 2\sum_{j=1}^\infty \lambda_j^2(\gamma),\\
\kappa_3(\gamma) & = & \displaystyle 8\sum_{j=1}^\infty \lambda_j^3(\gamma), & &
\kappa_4(\gamma) & = & \displaystyle 24\sum_{j=1}^\infty \lambda_j^3(\gamma),
\end{array}
$$
where the set $\{\lambda_j(\gamma)\}_{j\geq 1}$ denotes the solutions in $\lambda$ of equation \eqref{int:eq} with $\mu({\rm d}t)=w_\gamma(t){\rm d}t$.
To verify the accuracy of the solutions of equation \eqref{eq:det}, we calculated the approximations
\begin{equation} \label{approx:kum}
\begin{array}{rclcrcl}
\widehat{\kappa}_{n1}(\gamma) & = & \displaystyle \sum_{j=1}^m \widehat{\lambda}_j, & &
\widehat{\kappa}_{n2}(\gamma) & = & \displaystyle 2\sum_{j=1}^m \widehat{\lambda}_j^2,\\
\widehat{\kappa}_{n3}(\gamma) & = & \displaystyle 8\sum_{j=1}^m \widehat{\lambda}_j^3, & &
\widehat{\kappa}_{n4}(\gamma) & = & \displaystyle 8\sum_{j=1}^m \widehat{\lambda}_j^4,
\end{array}
\end{equation}
where $m$ is the number of non-null  solutions of equation  \eqref{eq:det} obtained by considering the first $n$ elements of the basis \eqref{eq:Hermite}.
Table \ref{tab:Ebner2021} presents $\kappa_i(\gamma)$ and $\widehat{\kappa}_{ni}(\gamma)$ for $\gamma=0.5, 1,2$, and $n=10(5)30$.


    \begin{table}[t] \centering
\begin{tabular}{cllllll}
$\gamma$ & & $n$ & $i=1$ & $i=2$ & $i=3$ & $i=4$\\ \hline
0.5 & $\kappa_i(\gamma)$          &    & 2.60125 & 4.71530 & 20.04364 & 133.19802     \\
    & $\widehat{\kappa}_{ni}(\gamma)$ & 10 & 2.57020 & 4.69630 & 19.97998 & 132.77912\\
    &                             & 15 & 2.59684 & 4.71475 & 20.04278 & 133.19340 \\
    &                             & 20 & 2.60065 & 4.71527 & 20.04361 & 133.19786 \\
    &                             & 25 & 2.60117 & 4.71530 & 20.04364 & 133.19807\\
    &                             & 30 & 2.60124 & 4.71530 & 20.04365 & 133.19809 \\ \hline
 1  & $\kappa_i(\gamma)$          &    & 7.78710e-1 & 5.43015e-1 & 8.71502e-1 & 2.20625 \\
    & $\widehat{\kappa}_{ni}(\gamma)$ & 10 & 7.77820e-1 & 5.42902e-1 & 8.71372e-1 & 2.20598\\
    &                             & 15 & 7.78679e-1 & 5.43014e-1 & 8.71502e-1 & 2.20625 \\
    &                             & 20 & 7.78709e-1 & 5.43015e-1 & 8.71502e-1 & 2.20625  \\
    &                             & 25 & 7.78710e-1 & 5.43015e-1 & 8.71502e-1 & 2.20625 \\
    &                             & 30 & 7.78710e-1 & 5.43015e-1 & 8.71502e-1 & 2.20625 \\ \hline
 2  & $\kappa_i(\gamma)$          &    & 2.02207e-1 & 4.57767e-2 & 2.39200e-2 & 1.97036e-2\\
    & $\widehat{\kappa}_{ni}(\gamma)$ & 10 & 2.02200e-1 & 4.57766e-2 & 2.39199e-2 & 1.97036e-2  \\
    &                             & 15 & 2.02207e-1 & 4.57767e-2 & 2.39200e-2 & 1.97036e-2 \\
    &                             & 20 & 2.02207e-1 & 4.57767e-2 & 2.39200e-2 & 1.97036e-2 \\
    &                             & 25 & 2.02207e-1 & 4.57767e-2 & 2.39200e-2 & 1.97036e-2\\
    &                             & 30 & 2.02207e-1 & 4.57767e-2 & 2.39200e-2 & 1.97036e-2
\end{tabular} \medskip
  \caption{First four cumulants and their Rayleigh-Ritz estimators of the distribution of $T(\gamma)=\|Z\|^2_{L^2}$, $Z$ a centered Gaussian process with covariance kernel $K_Z$ in \eqref{K_Z}.}
    \label{tab:Ebner2021}
\end{table}
    \item In \cite{HJM:2019}, Theorem 5.1, the covariance kernel is characterized as
    \begin{equation} \label{ET}
        C(s,t)=\exp(st)+\frac12(\exp(st)+\exp(-st))+2\cos(st)-st-4,\quad s,t\in\R.
    \end{equation}
    Note that the authors in this case assume $\gamma>1$.
     Table \ref{Hermite2} presents the estimators of the two largest eigenvalues of kernel \eqref{ET}, for $5\leq n \leq 13$, and $\gamma=1.5,2,3$. In this case, the convergence is notably quick, especially for $\gamma=3$.

\begin{table}[t]
    \centering
    \begin{tabular}{rrrrrrrrr}

   & \multicolumn{2}{c}{$\gamma=1.5$} &  & \multicolumn{2}{c}{$\gamma=2$} &  & \multicolumn{2}{c}{$\gamma=3$}\\
   \cline{2-3} \cline{5-6} \cline{8-9}
n &  \multicolumn{1}{c}{$\widehat{\lambda}_1$} &  \multicolumn{1}{c}{$\widehat{\lambda}_2$} &  &
      \multicolumn{1}{c}{$\widehat{\lambda}_1$}  &  \multicolumn{1}{c}{$\widehat{\lambda}_2$}
      &  &
      \multicolumn{1}{c}{$\widehat{\lambda}_1$}  &  \multicolumn{1}{c}{$\widehat{\lambda}_2$}\\ \hline
 5 & 3.367700e-1 & 1.860092e-1 & & 8.767026e-2 & 5.795029e-2 & & 1.388053e-2 & 1.271398e-2\\
 6 & 3.378871e-1 & 1.860092e-1 & & 8.769646e-2 & 5.795029e-2 & & 1.388069e-2 & 1.271398e-2\\
 7 & 3.378871e-1 & 1.861737e-1 & & 8.769646e-2 & 5.795333e-2 & & 1.388069e-2 & 1.271400e-2\\
 8 & 3.380118e-1 & 1.861737e-1 & & 8.769861e-2 & 5.795333e-2 & & 1.388070e-2 & 1.271400e-2\\
 9 & 3.380118e-1 & 1.861792e-1 & & 8.769861e-2 & 5.795335e-2 & & 1.388070e-2 & 1.271400e-2\\
10 & 3.380137e-1 & 1.861792e-1 & & 8.769862e-2 & 5.795335e-2 & & 1.388070e-2 & 1.271400e-2\\
11 & 3.380137e-1 & 1.861794e-1 & & 8.769862e-2 & 5.795335e-2 & & 1.388070e-2 & 1.271400e-2\\
12 & 3.380138e-1 & 1.861794e-1 & & 8.769862e-2 & 5.795335e-2 & & 1.388070e-2 & 1.271400e-2\\
13 & 3.380138e-1 & 1.861794e-1 & & 8.769862e-2 & 5.795335e-2 & & 1.388070e-2 & 1.271400e-2
 \end{tabular} \medskip
    \caption{Estimators of the two largest eigenvalues of kernel \eqref{ET} obtained for $5 \leq n \leq 13$ and $\gamma=1.5,2,3$  with the Rayleigh-Ritz method.}
    \label{Hermite2}
\end{table}

    \item In \cite{DEH:2020}, Theorem 5, the covariance kernel is characterized as
    \begin{eqnarray}\label{K_Norm}
        K(s,t)&=&\exp\left(-\frac{(s-t)^2}2\right)\left(
 \left(  \left( s-t \right) ^{2}-3 \right) ^{2}-6 \right) \\\nonumber&&+\exp\left(-\frac{s^2+t^2}2\right) \left( -\frac{{s}^{2}{t}^{2}
 \left( {s}^{2}-5 \right)  \left( {t}^{2}-5 \right)}2 +6({s}^{2}+{t}
^{2})-{s}^{4}-{t}^{4}\right.\\ && \nonumber\left.-{s}^{2}{t}^{2}-st \left( {s}^{2}-3 \right)
 \left( {t}^{2}-3 \right) -3 \right),
    \end{eqnarray}
    $s,t\in\R$. For $n=15$, the largest eigenvalue $\widehat{\lambda}_1(\gamma)$ is approximated by $\widehat{\lambda}_1(1/2)=0.966103408152626$, $\widehat{\lambda}_1(1)=0.600668091541773$, $\widehat{\lambda}_1(2)=0.398015644894253$, and
    $\widehat{\lambda}_1(3)=0.310469662783734$.

\end{itemize}

\end{example}

\subsection{Support $M = \mathbb{N}_0$}
If the domain of integration is $M=\mathbb{N}_0$, we concentrate on the measure with derivative with respect to the counting measure given by $w_\varrho(t)=e^{-\varrho}\varrho^t/t!$ for a positive parameter $\varrho>0$. In other words, we consider $\ell^2(\varrho)$, the (separable) Hilbert space of all infinite sequences $a=(a_0,a_1, \ldots)$ of complex numbers such that  $\sum_{t \geq 0}|a_t|^2w_\varrho(t)<\infty$, with the inner product defined as
\[
\langle a,b\rangle=\sum_{t \geq 0} a_t\bar{b}_tw_\varrho(t),
\]
for $a=(a_0,a_1, \ldots), \, b=(b_0,b_1, \ldots) \in \ell^2(\varrho)$. The practical computation of this inner product involves truncation of the infinite sum, say from $t=0$ up to $t=v$, for some positive integer $v$.

A suitable set of orthonormal polynomials can be expressed as
\begin{equation*}
\phi_k(x;\gamma)=\left(\varrho^k\,k!\right)^{1/2}
C_k(x; \varrho),
 \quad x\in\mathbb{N}_0,
\end{equation*}
where $C_k(x; \varrho)=
\sum_{\nu=0}^k(-1)^{k-\nu}\binom{k}{\nu}\nu!\varrho^{-\nu}\binom{x}{\nu}$ represents the Charlier polynomial of degree $k$.

\begin{example}(Testing for the von Mises distribution) \label{exm: von.Mises}
\begin{itemize}
    \item The covariance kernel in Theorem 1 of \cite{circular} is
    $$K(s,t)=\mathbb{E}\{\Upsilon(s,X;\theta) \overline{\Upsilon}(t,X;\theta)\}, \quad s,t \in \mathbb{N}_0,$$
    where
    \begin{eqnarray*}
    \Upsilon(s,X;\theta) & = &
    \cos(rX) -\Re \varphi (s;\theta) -\nabla \Re  \varphi (s;\theta)^\top L(X;\theta)\\
       &  & + {\rm i} \left\{\sin(rX) -\Im \varphi (s; \theta)-\nabla \Im \varphi (s;\theta)^\top L(X;\theta)\right\},
\end{eqnarray*}
$\theta$ is the vector of parameters of the law in the null hypothesis, $X$ is a circular random variable having the law in the null hypothesis with parameter vector $\theta$,
    $\Re \varphi (s;\theta)$ ($\Im \varphi (s;\theta)$) is the real (imaginary) part of the characteristic function of $X$,
   $\nabla \Re  \varphi (s;\theta)$ ($\nabla \Im \varphi (s;\theta)$) is the derivative of $\Re \varphi (s;\theta)$ ($\Im \varphi (s;\theta)$) with respect to $\theta$, $L(X;\theta)$ is the linear term in the Bahadur expansion of the estimator of $\theta$,
    and $\overline{\Upsilon}$ is the complex conjugate of ${\Upsilon}$.

    When applied to testing goodness-of-fit to the von Mises distribution, and the parameters $\mu$ (mean) and $\tau$ (concentration) are estimated with their maximum likelihood estimators, the above quantities become (see \cite{Mardia&Jupp2000})

    \begin{eqnarray*}
    \Re \varphi (s;\mu, \tau) & = &  \cos(\mu s) q(s; \tau),\\
    \Im \varphi (s;\mu, \tau) & = & \sin(\mu s) q(s, \tau),\\
    L(X;\mu, \tau) & = & \left(
    \begin{array}{c}
    \sin(x-\mu)/ q(1; \tau)\\  \{ \cos(x-\mu)-q(1; \tau)\}/\{1-q(1;\tau)^2-q(1;\tau)/\tau\} \end{array} \right),
    \end{eqnarray*}
    where $q(s; \tau)=I_s(\tau)/I_0(\tau)$  and $I_s(\tau)$ denotes the modified Bessel of the first kind and order $s$ (see, e.g. Chapter 9 of \cite{NIST}). After some calculations, one gets
   \begin{equation} \label{kernel.c}
   K(s,t) = \cos\left\{(s-t)\mu\right\} Q(s,t; \tau)+ {\rm i}  \sin\left\{(s-t)\mu\right\} Q(s,t; \tau), \quad s,t \in \mathbb{N}_0,
   \end{equation}
   where
    $$Q(s,t; \tau)  =  q(s-t; \tau)- q(s; \tau)q(t; \tau) \left\{1+st/(\tau q(1; \tau)\right\}
      -q'(s; \tau)q'(t, \tau)/\left\{1-q(1;\tau)^2-q(1;\tau)/\tau)\right\}$$
and $q'(s; \tau)=\frac{\partial}{\partial s} q(s; \tau)$.

Table \ref{tab:vM.RR} displays the estimators of the two largest eigenvalues obtained for $n=10,15,20$, $v=10$ (the same values were obtained for larger values of $v$) and $\mu=0$. As in the previous examples, rapid convergence is observed.

\begin{table}[t]
    \centering
    \begin{tabular}{ccccccccccc}
$\varrho$ &  $\tau$ & $n$ & $\widehat{\lambda}_1$ & $\widehat{\lambda}_2$  & & $\varrho$ &  $\tau$ & $n$ & $\widehat{\lambda}_1$ & $\widehat{\lambda}_2$\\
\hline
0.5 & 1 & 10 & 6.288772e-2 & 9.849891e-3 & & 0.5 & 5 & 10 & 8.438256e-3 & 5.328604e-3 \\
    &   & 15 & 6.288772e-2 & 9.849896e-3 & &     &   & 15 & 8.438256e-3 & 5.328604e-3\\
    &   & 20 & 6.288772e-2 & 9.849896e-3 & &     &   & 20 & 8.438256e-3 & 5.328604e-3\\ \hline
 1  & 1 & 10 & 1.632234e-1 & 4.686629e-2 & &  1  & 5 & 10 & 3.826515e-2 & 3.445026e-3\\
    &   & 15 & 1.632234e-1 & 4.686841e-2 & &     &   & 15 & 3.826515e-2 & 3.445029e-3\\
    &   & 20 & 1.632234e-1 & 4.686841e-2 & &     &   & 20 & 3.826515e-2 & 3.445029e-3
  \end{tabular} \medskip
    \caption{Estimators of the two largest eigenvalues  of kernel \eqref{kernel.c} obtained for $n=10,15,20$, $v=10$ and $\mu=0$,  with the Rayleigh-Ritz method. }
    \label{tab:vM.RR}
\end{table}

We juxtapose the results in Table \ref{tab:vM.RR} with those obtained by implementing the Monte Carlo procedure {outlined in Appendix \ref{app:MCA}}. The Monte Carlo procedure was repeated 500 times, for various values of $N$ (the number of generated samples).
Table  \ref{tab:vM.MC} showcases the mean and the standard deviation of the two largest values across the 500 replications. The means in Table \ref{tab:vM.MC}  approximate the values in  Table \ref{tab:vM.RR}, and they converge as $N$ increases. The standard deviations in Table \ref{tab:vM.MC}  diminish as $N$ increases.
From a computational perspective, the Rayleigh-Ritz method outperforms the Monte Carlo procedure, as the former computes the eigenvalues of a significantly smaller matrix.

\begin{table}[t]
    \centering
    \resizebox{\linewidth}{!} {
    \begin{tabular}{ccccrcrcccccrcr}
$\varrho$ &  $\tau$ & $N$ & \multicolumn{2}{c}{$\widehat{\lambda}_1$} & \multicolumn{2}{c}{$\widehat{\lambda}_2$}  & & $\varrho$ &  $\tau$ & $N$ &
\multicolumn{2}{c}{$\widehat{\lambda}_1$} & \multicolumn{2}{c}{$\widehat{\lambda}_2$}\\
\hline
0.5 & 1 &   50 &  6.349153e-2 & (26.8e-3) &  8.824334e-3 & (10.2e-3) & & 0.5 & 5 &   50 & 9.168138e-3 & (6.8e-3) & 2.903251e-4 & (5.3e-4)\\
    &   &  100 &  6.277916e-2 & (19.0e-3) &  9.430931e-3 &  (7.3e-3) & &     &   &  100 & 8.828349e-3 & (4.9e-3) & 3.828378e-4 & (5.2e-4)\\
    &   &  250 &  6.313665e-2 & (12.7e-3) &  9.840565e-3 &  (4.6e-3) & &     &   &  250 & 8.721612e-3 & (2.9e-3) & 4.852314e-4 & (4.3e-4)\\
    &   & 1000 &  6.286598e-2 &  (6.4e-3) &  9.850643e-3 &  (2.5e-3) & &     &   & 1000 & 8.483449e-3 & (1.4e-3) & 5.301540e-4 & (2.2e-4)\\
    &   & 2000 &  6.269159e-2 &  (4.5e-3) &  9.840051e-3 &  (1.8e-3) & &     &   & 2000 & 8.538285e-3 & (1.0e-3) & 5.430429e-4 & (1.6e-4)\\
    &   & 3000 &  6.270766e-2 &  (3.7e-3) &  9.818919e-3 &  (1.5e-3) & &     &   & 3000 & 8.487909e-3 & (1.0e-3) & 5.405830e-4 & (1.6e-4)\\
    &   & 4000 &  6.282082e-2 &  (3.6e-3) &  9.825737e-3 &  (1.3e-3) & &     &   & 4000 & 8.471089e-3 & (0.7e-3) & 5.362364e-4 & (1.1e-4) \\
    &   & 5000 &  6.288797e-2 &  (2.8e-3) &  9.787226e-3 &  (1.1e-3) & &     &   & 5000 & 8.453944e-3 & (0.6e-3) & 5.351108e-4 &  (1-0e-4)\\
\hline
1 & 1 &   50 & 1.663794e-1 & (43.0e-3) & 4.383321e-2 & (18.9e-3) & & 1 & 5 &   50 & 3.826183e-2 & (13.9e-3) & 2.835939e-3 & (26.6e-4) \\
  &   &  100 & 1.648465e-1 & (28.6e-3) & 4.481114e-2 & (13.6e-3) & &   &   &  100 & 3.836605e-2 & (10.0e-3) & 3.294036e-3 & (21.8e-4) \\
  &   &  250 & 1.631288e-1 & (17.6e-3) & 4.576026e-2 &  (8.7e-3) & &   &   &  250 & 3.830483e-2 &  (6.7e-3) & 3.339658e-3 & (13.6e-4) \\
  &   & 1000 & 1.631152e-1 &  (9.3e-3) & 4.652113e-2 &  (4.4e-3) & &   &   & 1000 & 3.813966e-2 &  (3.3e-3) & 3.434499e-3 &  (6.4e-3) \\
  &   & 2000 & 1.629177e-1 &  (6.4e-3) & 4.654794e-2 &  (3.2e-3) & &   &   & 2000 & 3.826935e-2 &  (2.4e-3) & 3.452576e-3 &  (4.6e-3) \\
  &   & 3000 & 1.631266e-1 &  (5.5e-3) & 4.661693e-2 &  (2.6e-3) & &   &   & 3000 & 3.832378e-2 &  (1.9e-3) & 3.464659e-3 &  (3.9e-4) \\
  &   & 4000 & 1.631981e-1 &  (4.8e-3) & 4.674310e-2 &  (2.1e-3) & &   &   & 4000 & 3.834800e-2 &  (1.7e-3) & 3.462685e-3 &  (3.4e-4)  \\
  &   & 5000 & 1.632090e-1 &  (4.3e-3) & 4.678101e-2 &  (1.9e-3) & &   &   & 5000 & 3.832935e-2 &  (1.5e-3) & 3.459746e-3 &  (3.1e-4)
   \end{tabular}} \medskip
    \caption{Mean and standard deviation (in parenthesis) in 500 replications of the Monte Carlo method for estimating the  two largest eigenvalues   of kernel \eqref{kernel.c},  for several values of $N$.}
    \label{tab:vM.MC}
\end{table}

\end{itemize}
\end{example}


\subsection{Support $M=\mathbb{R}^d$}
In this scenario, we concentrate on the weight function $w_\gamma(t)=\exp(-\gamma\|t\|^2)$, $t\in\R^d$, $\gamma>0$, where $\|x\|=\sqrt{x^\top x}$ represents the Euclidean norm and $\top$ denotes the transpose of a vector. In the spirit of \cite{B:1996}, a suitable set of orthonormal polynomials can be expressed using \eqref{eq:Hermite} as
\begin{equation} \label{basis_d}
    \widetilde{\phi}_{k_1,\ldots, k_d}(x;\gamma)=\prod_{j=1}^d\phi_{k_j}(x_j)=\prod_{j=1}^d\left(2^{k_j}\,k_j!\sqrt{\pi/\gamma}\right)^{-d/2}H_{k_j}\left(\sqrt{\gamma}x_j\right),\quad x=(x_1,\ldots,x_d)\in\R^d, \quad k_1, \ldots, k_d \in \mathbb{N}_0.
\end{equation}
We focus our attention on the classic BHEP test as presented in \cite{HW:1997} with covariance kernel given in display (2.3)
\begin{equation} \label{kernel:BHEP}
    K(s,t)=\exp\left(-\frac{\|s-t\|^2}{2}\right)-\left(1+s^\top t+\frac{(s^\top t)^2}2\right)\exp\left(-\frac{\|s\|^2+\|t\|^2}{2} \right),\quad s,t\in\R^d.
\end{equation}
The weight function considered in  \cite{HW:1997} is
\[
w(t)=\varphi_\beta(t)=(2\pi\beta^2)^{-d/2} \exp \left( -\frac{ \|t\|^2}{2 \beta^2}\right), \quad t\in\mathbb{R}^d, \quad \beta>0.
\]
Clearly, for $\gamma=1/(2  \beta^2)$, we have that
\begin{equation} \label{eq1}
w_\gamma(t)=(2 \pi \beta^2)^{d/2} \varphi_\beta(t),
\end{equation}
and in such a case,  from \eqref{eq1},  
$$ \int K(s,t) f(s) \varphi_\beta(s) ds= \theta f(t)$$
if and only if
$$ \int K(s,t) f(s) w_\gamma(s) ds= \lambda f(t), \quad \mbox{with} \quad \lambda=(2 \pi \beta^2)^{d/2} \theta.$$
The exact expression of the first three cumulants of the law of the distribution of $T(\beta)=\|Z\|^2_{L^2}$ with $w=\varphi_\beta$ (denote them as $\kappa_1(\beta)$, $\kappa_2(\beta)$ and $\kappa_3(\beta)$)  have been given in  \cite{HW:1997}. If such cumulants are denoted as $\kappa_1(\gamma)$, $\kappa_2(\gamma)$ and $\kappa_3(\gamma)$, respectively, when $w=w_\gamma$, then from \eqref{eq1}, for $\gamma=1/(2  \beta^2)$, we have that
$$ \kappa_1(\gamma)  =  (\pi/\gamma)^{d/2} \kappa_1(\beta), \quad \kappa_2(\gamma)  =  (\pi/\gamma)^{d} \kappa_2(\beta), \quad  \kappa_3(\gamma)  =  (\pi/\gamma)^{3d/2} \kappa_3(\beta).$$

In order to check the accuracy of the solutions of equation \eqref{eq:det}, we calculated the exact values of  $\kappa_1(\gamma)$, $\kappa_2(\gamma)$ and $\kappa_3(\gamma)$ and the approximations in
\eqref{approx:kum}, where now $m$ is the number of non-null  solutions of equation  \eqref{eq:det} obtained by considering all elements in the basis \eqref{basis_d} with $k_1+\ldots+k_d \leq n$. Tables \ref{tab:dim1}--\ref{tab:dim3} display $\kappa_i(\gamma)$ and $\widehat{\kappa}_{ni}(\gamma)$ for $\gamma=0.5, 1,2$,  $n=10(5)30$ and $d=1,2,3$. Looking at these tables it can be concluded that $\kappa_i(\gamma)$ and $\widehat{\kappa}_{ni}(\gamma)$ are very close for $n\geq 15$, specially when $\gamma=1,2$.

Note that using the formulas from the first point in Example \ref{exm: normality} we can easily approximate the hitherto unknown fourth cumulant of the limit distribution of the BHEP test. This would allow to fit a Pearson system of distributions as described in Section \ref{Pearson.System} to efficiently approximate the critical values of the test statistic.
\begin{table} \centering
\begin{tabular}{clllll}
$\gamma$ & & $n$ & $i=1$ & $i=2$ & $i=3$\\ \hline
0.5 & $\kappa_i(\gamma)$          & true   & 3.35825e-1 & 9.57324e-2 & 6.30525e-2 \\
    & $\widehat{\kappa}_{ni}(\gamma)$ & 10 & 3.30722e-1 & 9.41680e-2 & 6.16292e-2\\
    &                             & 15 & 3.35124e-1 & 9.56358e-2 & 6.29862e-2 \\
    &                             & 20 & 3.35728e-1 & 9.57275e-2 & 6.30492e-2  \\
    &                             & 25 & 3.35813e-1 & 9.57322e-2 & 6.30524e-2 \\
    &                             & 30 & 3.35823e-1 & 9.57324e-2 & 6.30525e-2\\ \hline
 1  & $\kappa_i(\gamma)$          &  true  & 8.83130e-2 & 8.07001e-3 & 1.69471e-3 \\
    & $\widehat{\kappa}_{ni}(\gamma)$ & 10 & 8.81463e-2 & 8.05925e-3 & 1.69165e-3\\
    &                             & 15 & 8.83081e-2 & 8.06992e-3 & 1.69469e-3 \\
    &                             & 20 & 8.83128e-2 & 8.07001e-3 & 1.69471e-3  \\
    &                             & 25 & 8.83130e-2 & 8.07001e-3 & 1.69471e-3 \\
    &                             & 30 & 8.83130e-2 & 8.07001e-3 & 1.69471e-3 \\ \hline
 2  & $\kappa_i(\gamma)$          &  true  & 1.67944 & 3.51349e-4 & 1.68584e-5  \\
    & $\widehat{\kappa}_{ni}(\gamma)$ & 10 & 1.67930 & 3.51336e-4 & 1.68575e-5 \\
    &                             & 15 & 1.67944 & 3.51349e-4 & 1.68584e-5 \\
    &                             & 20 & 1.67944 & 3.51349e-4 & 1.68584e-5 \\
    &                             & 25 & 1.67944 & 3.51349e-4 & 1.68584e-5\\
    &                             & 30 & 1.67944 & 3.51349e-4 & 1.68584e-5
\end{tabular} \medskip
  \caption{
  First three cumulants and their Rayleigh-Ritz estimators of the distribution of $T(\gamma)=\|Z\|^2_{L^2}$, $Z$ a centered Gaussian process with covariance kernel $K$ in \eqref{kernel:BHEP} and $d=1$.}
    \label{tab:dim1}
\end{table}

\begin{table} \centering
\begin{tabular}{clllll}
$\gamma$ & & $n$ & $i=1$ & $i=2$ & $i=3$ \\ \hline
0.5 & $\kappa_i(\gamma)$          & true   & 1.86169 & 7.38190e-1 & 7.87273e-1\\
    & $\widehat{\kappa}_{ni}(\gamma)$ & 10 & 1.80572 & 7.24167e-1 & 7.67663e-1\\
    &                             & 15 & 1.85336 & 7.37210e-1 & 7.86174e-1\\
    &                             & 20 & 1.86047 & 7.38133e-1 & 7.87209e-1 \\
    &                             & 25 & 1.86152 & 7.38187e-1 & 7.87270e-1\\
    &                             & 30 & 1.86166 & 7.38190e-1 & 7.87272e-1\\ \hline
 1  & $\kappa_i(\gamma)$          & true   & 3.92699e-1 & 4.36479e-2 & 1.28019e-2 \\
    & $\widehat{\kappa}_{ni}(\gamma)$ & 10 & 3.91340e-1 & 4.35734e-2 & 1.27726e-2\\
    &                             & 15 & 3.92654e-1 & 4.36472e-2 & 1.28018e-2\\
    &                             & 20 & 3.92698e-1 & 4.36479e-2 & 1.28019e-2\\
    &                             & 25 & 3.92699e-1 & 4.36479e-2 & 1.28019e-2\\
    &                             & 30 & 3.92699e-1 & 4.36479e-2 & 1.28019e-2\\ \hline
 2  & $\kappa_i(\gamma)$          & true   & 5.81776e-2 & 1.21539e-3 & 6.65175e-5\\
    & $\widehat{\kappa}_{ni}(\gamma)$ & 10 & 5.81692e-2 & 1.21533e-3 & 6.65130e-5\\
    &                             & 15 & 5.81776e-2 & 1.21539e-3 & 6.65175e-5\\
    &                             & 20 & 5.81776e-2 & 1.21539e-3 & 6.65175e-5\\
    &                             & 25 & 5.81776e-2 & 1.21539e-3 & 6.65175e-5\\
    &                             & 30 & 5.81776e-2 & 1.21539e-3 & 6.65175e-5
\end{tabular} \medskip
  \caption{
  First three cumulants and their Rayleigh-Ritz estimators of the distribution of $T(\gamma)=\|Z\|^2_{L^2}$, $Z$ a centered Gaussian process with covariance kernel $K$ in \eqref{kernel:BHEP} and $d=2$.}
    \label{tab:dim2}
\end{table}

\begin{table} \centering
\begin{tabular}{clllll}
$\gamma$ & & $n$ & $i=1$ & $i=2$ & $i=3$\\ \hline
0.5 & $\kappa_i(\gamma)$          & true   & 7.16174 & 3.89456 & 6.64335 \\
    & $\widehat{\kappa}_{ni}(\gamma)$ & 10 & 6.78990 & 3.80774 & 6.46828\\
    &                             & 15 & 7.09987 & 3.88788 & 6.63171\\
    &                             & 20 & 7.15194 & 3.89413 & 6.64257 \\
    &                             & 25 & 7.16030 & 3.89453 & 6.64331\\
    &                             & 30 & 7.16153 & 3.89456 & 6.64335\\ \hline
 1  & $\kappa_i(\gamma)$          & true   & 1.20027 & 1.61360e-1 & 6.53432e-2\\
    & $\widehat{\kappa}_{ni}(\gamma)$ & 10 & 1.19318 & 1.61016e-1 & 6.51617e-2\\
    &                             & 15 & 1.20027 & 1.61356e-1 & 6.53418e-2\\
    &                             & 20 & 1.20026 & 1.61360e-1 & 6.53432e-2 \\
    &                             & 25 & 1.20027 & 1.61360e-1 & 6.53432e-2\\
    &                             & 30 & 1.20027 & 1.61360e-1 & 6.53432e-2\\ \hline
 2  & $\kappa_i(\gamma)$          & true   & 1.38008e-1 & 2.87382e-3  & 1.77339e-4 \\
    & $\widehat{\kappa}_{ni}(\gamma)$ & 10 & 1.37974e-1 & 2.87363e-3  & 1.77323e-4\\
    &                             & 15 & 1.38008e-1 & 2.87382e-3  & 1.77339e-4\\
    &                             & 20 & 1.38008e-1 & 2.87382e-3  & 1.77339e-4\\
    &                             & 25 & 1.38008e-1 & 2.87382e-3  & 1.77339e-4\\
    &                             & 30 & 1.38008e-1 & 2.87382e-3  & 1.77339e-4
\end{tabular} \medskip
  \caption{First three cumulants and their Rayleigh-Ritz estimators of the distribution of $T(\gamma)=\|Z\|^2_{L^2}$, $Z$ a centered Gaussian process with covariance kernel $K$ in \eqref{kernel:BHEP} and $d=3$.
}
    \label{tab:dim3}
\end{table}

\section{Usefulness of eigenvalues estimation}\label{sec:Con}
In this section, we introduce two statistical use cases where an accurate approximation of the eigenvalues of the covariance operator is of paramount importance.

\subsection{Distribution approximation}
\label{Pearson.System}

As noted in the Introduction, if $Z$ is a Gaussian process that takes values in an appropriate function space and $\|\cdot\|$ is a corresponding norm, then $\|Z\|^2$ follows the distribution $W=\sum_{j=1}^\infty \lambda_j N_j^2$, where $N_1, N_2, \ldots$ are iid standard normal random variables, and $\{\lambda_j\}_{j\in\N}$ is the sequence of positive eigenvalues of the covariance operator, ${\cal K}$, of $Z$. Therefore, if $\widehat{\lambda}_1, \ldots, \widehat{\lambda}_m$ are estimators of the $m$ largest eigenvalues of ${\cal K}$, the distribution of $\|Z\|^2$ can be approximated by that of $\widehat{W}_m=\sum_{j=1}^m \widehat{\lambda}_j N_j^2$. The distribution of a linear combination of chi-squared variates can subsequently be approximated by simulation or by employing the Imhof method \cite{Imhof}, which is implemented in the {\tt R} package {\tt CompQuadForm} \cite{CompQuadForm}. This procedure has been utilized in \cite{NJ:2016,  RJ:2021, RJ:2019} with eigenvalue estimators obtained as described in Appendix \ref{app:MCA}.

Next, we observe that by approximating the eigenvalues of the covariance operator, we can approximate the limiting distribution by a Pearson system. Initially, by directly evaluating integrals, the first two cumulants of the distribution of $\|Z\|_{L^2}^2$ are
\begin{equation*}
\kappa_{1}=\E \|Z\|_{L^2}^2=\int_M K(t,t)\,w(t)\,\mbox{d}t
\end{equation*}
and
\begin{equation*}
\kappa_{2}=\mathbb{V}(\|Z\|_{L^2}^2)=2\int_M\int_M K^2(s,t)\,w(s)w(t)\,\mbox{d}s\mbox{d}t,
\end{equation*}
where $\mathbb{V}(\cdot)$ represents the variance. Following the methodology in \cite{H:1990,S:1976}, the third and fourth cumulants can be computed as
\begin{equation*}
\kappa_{j}=2^{j-1}(j-1)!\int_M K_j(t,t) w(t)\,\mbox{d}t,
\end{equation*}
where $K_j(s,t)$, the $j^{\text{th}}$ iterate of $K(s,t)$, is defined as
\begin{eqnarray*}
K_{j}(s,t)&=&\int_M K_{j-1}(s,u)K_0(u,t)w(u)\,\mbox{d}u,\quad j\ge2,\\
K_1(s,t)&=&K(s,t).
\end{eqnarray*}
Utilizing the identities in Remark 2.2 of \cite{E:2023} (also compare to Corollary 1.3 in \cite{DM:2003}), we observe that
\begin{equation*}
\sum_{k=1}^\infty\lambda^j_k=\kappa_{j}/(2^{j-1}(j-1)!),\quad j=1,2,3,4,
\end{equation*}
thus sums of powers of the eigenvalues approximate the cumulants of the limiting distribution. This approach, applied by direct calculation of $\kappa_1, \kappa_2, \kappa_3,\kappa_4$ as described above, has been utilized in \cite{ E:2020, EH:2023, H:1990}.

\subsection{Bahadur efficiency}

Here we demonstrate the importance of the approximation of eigenvalues for the test's quality assessment, by calculating local Bahadur efficiency - one of  often-used asymptotic criteria.
Its widespread adoption as an asymptotic quality measure is primarily due to its applicability to test statistics with non-normal limiting distributions.
For details, we refer to \cite{nikitinKnjiga}, while a brief review is given in Appendix \ref{ap: LABE}.


Let $\mathcal{G}=\{G_{\theta}(x),\;\theta>0\}$ be a family of alternative distribution functions, with $G_0(x)$ being the null family of distributions.
The relative Bahadur efficiency of two sequences of test statistics $\{T_n\}$ and $\{V_n\}$ can be represented as the ratio of their
approximate Bahadur slope functions $c^*_T(\theta)$ and $c^*_V(\theta)$, which are associated with the rate of exponential decrease for
the level of significance achieved under the alternative. Typically, we are interested in the case when $\theta\to 0$, i.e., when we are comparing the behavior of tests against nearby alternatives.

Given that likelihood ratio (LR) tests are optimal in the Bahadur sense (refer to \cite{bahadur1967optimal, rublik1989optimality}), for close alternatives from $\mathcal{G}$, the absolute local approximate
Bahadur efficiency for $T_n$ is defined as the ratio
of $c_T(\theta)$ and the corresponding
slope of the LR test, which equals $2K (\theta)$ -- twice the Kullback–Leibler (KL) distance from the given alternative to the class of distributions within the null hypothesis (see also \cite{meintanis2022bahadur}).

If the weak limit of $T_n$ is $||Z||^2_{L^2}$, under certain regularity conditions (see \cite{meintanis2022bahadur}) the local approximate Bahadur slope is equal to $c^\ast_T(\theta)={b^{''}_T(0)}\theta^2/({2\lambda_1})+o(\theta^2),\;\theta\to 0.$ The coefficient $b_T(\theta)$ is the limit in probability of ${T_n}/{n}$. Therefore, the computation of $\lambda_1$ is crucial in the derivation of local approximate Bahadur efficiency.



\begin{example}
In this example, we bridge the gap in the literature by presenting Bahadur efficiencies of exponentiality test associated with kernel $\rho$ in \eqref{kernel rho}, against the following frequently considered close alternatives:

\begin{itemize}
		\item the Weibull distribution with density
		\begin{equation*}
		 g_\theta(x)=e^{-x^{1+\theta}}(1+\theta)x^\theta,\quad \theta>0,\quad x\geq0;
		\end{equation*}
		\item the gamma distribution with density
		\begin{equation*}
		g_\theta(x)=\frac{x^\theta e^{-x}}{\Gamma(\theta+1)}, \quad  \theta>0 \quad x\geq0;
		\end{equation*}
		\item
		the Makeham distribution with density
		\begin{equation*}
		g_\theta(x)=e^{-x-\theta e^{x}}(1+\theta e^x),\quad \theta>0, \quad x\geq0;
		\end{equation*}
		\item the linear failure rate (LFR) distribution with density
		\begin{equation*}
		g_\theta(x)=e^{-x-\theta\frac{x^2}{2}}(1+\theta x), \quad \theta>0 \quad ,x\geq0;
		\end{equation*}
		\item the mixture of exponential distributions with negative weights (EMNW($\beta$)) with density
		\begin{equation*}
		g_\theta(x;\beta)=(1+\theta)e^{-x}-\theta\beta e^{-\beta x},\quad  \theta\in\left(0,\frac{1}{\beta-1}\right], \quad x\geq0.
		\end{equation*}
	\end{itemize}
	
Efficiencies of the test associated with kernel $\rho$ 
are already presented in \cite{E:2023},
while a comprehensive comparison of exponentiality tests in terms of local approximate Bahadur efficiencies can be found in \cite{cuparic2022}. For some recent results, see also \cite{meintanis2022bahadur}. Therefore, as a complement of previous results, in Table \ref{tab:BE exp} we just present results for the test from \cite{K:2001}. The Bahadur efficiencies are obtained using the largest eigenvalues obtained via Raylegh-Ritz method for $n=10$.


\begin{table*}[t]
    \centering
    \begin{tabular}{lc|cccc}
        &  &  \multicolumn{4}{c}{$T_{\gamma}$ with kernel \eqref{K2.Laguerre}}\\
         Alternative & $\gamma$ & $0$ & 1 & 2 & 3 \\\hline
         Weibull&&0.672&0.817&0.868&0.884\\
         Gamma&&0.453&0.625&0.722&0.781\\
         Makeham&&0.855&0.987&0.982&0.835\\
         LFR&&0.962&0.798&0.657&0.555\\
         EMNW(3)&&0.668&0.889&0.974&0.995
    \end{tabular} \medskip

    \caption{BE of exponentially of tests from \cite{K:2001}
    \label{tab:BE exp} in Example \ref{exm: exp}}
\end{table*}

\end{example}

\begin{example}
Here we reconsider tests from Example \ref{exm: normality} and calculate local approximate Bahadur efficiencies against the following commonly used alternatives  (see \cite{milovsevic2021bahadur} and \cite{meintanis2022bahadur}). For the calculation of efficiencies, we use the largest eigenvalues obtained via Raylegh-Ritz method for $n=10$. The results are presented in Table \ref{eff}.

\begin{itemize}
\item the Lehmann alternatives
\begin{equation*}
        g^{(1)}_\theta(x)=(1+\theta)F_{\vartheta_0}^{\theta}(x)f_{\vartheta_0}(x)
    \end{equation*}
    \item the first Ley-Paindaveine alternatives
    \begin{equation*}
        g^{(2)}_\theta(x)=f_{\vartheta_0}(x) e^{-\theta(1-F_{\vartheta_0}(x))}
        \big(1+\theta F_{\vartheta_0}(x)\big)
    \end{equation*}
    \item the second Ley-Paindaveine alternatives
    \begin{equation*}
        g^{(3)}_\theta(x)=f_{\vartheta_0}(x)\big(1-\theta\pi\cos(\pi F_{\vartheta_0}(x)\big)
    \end{equation*}
    \item the contamination alternatives
    \begin{equation*}
        g^{(4)}_\theta(x;\mu,\sigma^2)=(1-\theta)f_{\vartheta_{0}}(x)+\theta \frac1{\sigma}f_{\vartheta_0}
        \Big(\frac{x-\mu}{\sigma}\Big) \label{contamination}
    \end{equation*}
\end{itemize}

The outcomes are displayed in Table \ref{eff}. It is evident that the tests are quite sensitive to the selection of the tuning parameter, and as a result, the choice of the weight function significantly influences the behaviour of the tests. The test proposed in \cite{DEH:2020} is the least affected by the choice of tuning parameter. This test also demonstrates exceptional efficiencies against Lehmann's alternative, which are very close to the highest efficiencies of the energy and BHEP tests (refer to \cite{meintanis2022bahadur}). However, this test is considerably less efficient under the contamination with $N(0,0.5)$ alternative. Conversely, in this scenario, despite being the least efficient against many alternatives, the test proposed in \cite{HJM:2019} has exhibited superior performance compared to all the tests considered in the literature.

\begin{table}
    \centering
    \begin{tabular}{lc|cccc|ccc|cccc}
      &  &  \multicolumn{4}{c}{$T_\gamma$ with kernel \eqref{K_Z}}&  \multicolumn{3}{c}{$T_\gamma$ with kernel  \eqref{ET}}  &   \multicolumn{4}{c}{$T_\gamma$ with
 kernel \eqref{K_Norm}}   \\
     Alternative & $\gamma$ & $1/2$ & 1 & 2 & 3  &$3/2$& 2 & 3 & $1/2$ & 1 & 2 & 3  \\ \hline
       $g^{(1)}$ &  &0.637  &0.809  & 0.918 & 0.951  &{ 0.386}  &{ 0.543}  &{ 0.831
       }  &0.916  &0.964  & 0.959&0.964 \\
       $g^{(2)}$ & &0.862  &0.976  &0.996  & 0.984    &0.256  & 0.387 &0.697  &0.908  &0.968  &0.980&0.975 \\
       $g^{(3)}$  &  &0.955  & 0.986 & 0.939 &0.899  &0.194  &0.299  &0.387  &0.802  &0.875  & 0.895& 0.883 \\
       $g^{(4)}(1,1)$  &  &0.419  & 0.555 &0.658  &0.697  & 0.564&0.698  &0.924  &0.722  & 0.752 & 0.731&0.730\\
       $g^{(4)}(0.5,1)$  &  &0.552  & 0.719 &0.838  &0.880  & {0.449}  &{0.604}  &{0.883}  & 0.862 &0.905  &0.895&0.900 \\
        $g^{(4)}(0,0.5)$  &  &0.552  &0.560  &0.362  &0.261  &{0.823}    &{0.823}  & {0.809} &0.550  & 0.502 &0.396&0.312 \\
    \end{tabular}
    \medskip
    \caption{BE of normality tests in Example \ref{exm: normality}}
    \label{eff}
\end{table}

\end{example}

\section{Comments and Outlook}

As mentioned in the Introduction, a primary motivation for approximating the eigenvalues of a kernel is related to goodness-of-fit testing. For location-scale families, the kernel associated with the test statistic is usually free from unknown parameters. However, for families with a shape parameter, the kernel typically depends on this parameter, which is often unknown in practice. If the kernel is continuous as a function of the parameter, the result in Theorem \ref{thm} remains valid if the parameter is replaced by a consistent estimator. The convergence is in probability or almost surely, depending on whether the parameter estimator is weakly or strongly consistent, respectively. Clearly, this implies dependence of the eigenvalues on the parameter.

We conclude this paper by emphasizing the high adaptability of the introduced method. It is applicable in any dimension, in discrete scenarios, and even for relatively small $n$, providing reasonable deterministic approximations. This is a clear advantage over the competing stochastic method outlined in Appendix \ref{app:AAA}. The Rayleigh-Ritz method can be extended to supports $M$ that model more complex data sets, such as matrix-valued data or functional data in a Hilbert space. In such cases, it is necessary to find an appropriate set of orthonormal basis functions with respect to the integrating measure $\mu$. Examples of weighted $L^2$-type statistics in these spaces and the corresponding kernels of Gaussian processes is found in \cite{HR:2020} and \cite{HJ:2021}. Fortunately, the method can be efficiently computed if the orthonormal basis for the underlying $L^2$ space can be implemented and works well even for complex kernels. While analytical solutions to the eigenvalue problem are always preferred, the complexity of the formulas provided for the first part of the kernel of the BHEP test (see \cite{EH:2023}) suggests that such solutions will be rare for more intricate kernels. Thus, a proficient approximation provided by the Rayleigh-Ritz method will suffice for the time being.


%


\color{black}
\section*{Acknowledgement}
The research of B. Ebner and B.  Milo\v{s}evi\'{c} was supported by the bilateral cooperation project ``Modeling complex data - Selection and Specification” between the Federal Republic of Germany and the Republic of Serbia (the contract 337-00-19/2023-01/6). The project underlying this article was funded by the Federal Ministry of Education and Research of Germany, the Ministry of Science, Technological Development and Innovation of Republic of Serbia (the contract 451-03-66/2024-03/200104), and COST action CA21163 - Text, functional and other high-dimensional data in econometrics: New models, methods, applications (HiTEc). The research of M.D. Jim\'{e}nez-Gamero was supported by grant PID2023-147058NA-I00 (Ministerio de Ciencia, Innovaci\'on y Universidades, Spain).

\bibliographystyle{abbrv}
\bibliography{lit_EVKO}

\begin{thebibliography}{10}

\bibitem{AES:2023}
J.~S. Allison, B.~Ebner, and M.~Smuts.
\newblock Logistic or not logistic?
\newblock {\em Statistica Neerlandica}, 77(4):429--443, 2023.

\bibitem{AP:2017}
J.~S. Allison and C.~Pretorius.
\newblock A monte carlo evaluation of the performance of two new tests for
  symmetry.
\newblock {\em Computational Statistics}, 32(4):1323--1338, 2017.

\bibitem{bahadur1967optimal}
R.~Bahadur.
\newblock An optimal property of the likelihood ratio statistic.
\newblock In {\em Proceedings of the Fifth Berkeley Symposium on Mathematical
  Statistics and Probability}, volume~1, pages 13--26, 1967.

\bibitem{BBQ:2018}
T.~Bahraoui, T.~Bouezmarni, and J.-F. Quessy.
\newblock A family of goodness-of-fit tests for copulas based on characteristic
  functions.
\newblock {\em Scandinavian Journal of Statistics}, 45(2):301--323, 2018.

\bibitem{B:1977}
C.~T.~H. Baker.
\newblock {\em The numerical treatment of integral equations}.
\newblock Oxford University Press, 1977.

\bibitem{B:1996}
L.~Baringhaus.
\newblock Fibonacci numbers, {L}ucas numbers and integrals of certain
  {G}aussian processes.
\newblock {\em Proceedings of the American Mathematical Society},
  124(12):3875--3884, 1996.

\bibitem{BEH:2017}
L.~Baringhaus, B.~Ebner, and N.~Henze.
\newblock The limit distribution of weighted ${L}^2$-goodness-of-fit statistics
  under fixed alternatives, with applications.
\newblock {\em Annals of the Institute of Statistical Mathematics},
  69(5):969--995, 2017.

\bibitem{BGT:2018}
L.~Baringhaus, D.~Gaigall, and J.~P. Thiele.
\newblock Statistical inference for {$L^2$}-distances to uniformity.
\newblock {\em Computational Statistics}, 33(4):1863--1896, 2018.

\bibitem{BH:2008}
L.~Baringhaus and N.~Henze.
\newblock A new weighted integral goodness-of-fit statistic for exponentiality.
\newblock {\em Statistics \& Probability Letters}, 78(8):1006--1016, 2008.

\bibitem{BE:2019b}
S.~Betsch and B.~Ebner.
\newblock A new characterization of the gamma distribution and associated
  goodness-of-fit tests.
\newblock {\em Metrika}, 82(7):779--806, 2019.

\bibitem{BL:2005}
M.~Bilodeau and P.~L. de~Micheaux.
\newblock A multivariate empirical characteristic function test of independence
  with normal marginals.
\newblock {\em Journal of Multivariate Analysis}, 95(2):345 -- 369, 2005.

\bibitem{bozin2020new}
V.~Bo{\v{z}}in, B.~Milo{\v{s}}evi{\'c}, Y.~Y. Nikitin, and M.~Obradovi{\'c}.
\newblock \color{black}{N}ew characterization-based symmetry tests.
\newblock {\em Bulletin of the Malaysian Mathematical Sciences Society},
  43(1):297--320, 2020.

\bibitem{cuparic2022}
M.~Cupari{\'c}, B.~Milo{\v{s}}evi{\'c}, and M.~Obradovi{\'c}.
\newblock New consistent exponentiality tests based on {V}-empirical {L}aplace
  transforms with comparison of efficiencies.
\newblock {\em Revista de la Real Academia de Ciencias Exactas, F{\'\i}sicas y
  Naturales. Serie A. Matem{\'a}ticas}, 116(1):1--26, 2022.

\bibitem{DM:2003}
P.~Deheuvels and G.~Martynov.
\newblock Karhunen-{L}o{\`e}ve expansions for weighted {W}iener processes and
  {B}rownian bridges via {B}essel functions.
\newblock In J.~Hoffmann-J{\o}rgensen, J.~A. Wellner, and M.~B. Marcus,
  editors, {\em High Dimensional Probability III}, pages 57--93, Basel, 2003.
  Birkh{\"a}user Basel.

\bibitem{DW:2014}
J.~Duan and W.~Wang.
\newblock {\em Effective Dynamics of Stochastic Partial Differential
  Equations}.
\newblock Academic Press, Boston, 2014.

\bibitem{CompQuadForm}
P.~Duchesne and P.~L. {de Micheaux}.
\newblock Computing the distribution of quadratic forms: Further comparisons
  between the {L}iu-{T}ang-{Z}hang approximation and exact methods.
\newblock {\em Computational Statistics and Data Analysis}, 54:858--862, 2010.

\bibitem{DSpartII}
N.~Dunford and J.~T. Schwartz.
\newblock {\em Linear operators. {P}art {II}: {S}pectral theory. {S}elf adjoint
  operators in {H}ilbert space}.
\newblock Wiley, New York-London, 1963.

\bibitem{DEH:2020}
P.~Dörr, B.~Ebner, and N.~Henze.
\newblock Testing multivariate normality by zeros of the harmonic oscillator in
  characteristic function spaces.
\newblock {\em Scandinavian Journal of Statistics}, 48(2):456--501, 2021.

\bibitem{E:2020}
B.~Ebner.
\newblock On combining the zero bias transform and the empirical characteristic
  function to test normality.
\newblock {\em ALEA}, 18:1029--1045, 2020.

\bibitem{E:2023}
B.~Ebner.
\newblock The test of exponentiality based on the mean residual life function
  revisited.
\newblock {\em Journal of Nonparametric Statistics}, 35(3):601--621, 2023.

\bibitem{EH:2023b}
B.~Ebner and N.~Henze.
\newblock Bahadur efficiencies of the {E}pps--{P}ulley test for normality.
\newblock {\em Journal of Mathematical Sciences}, 273(5):861--870, Jul 2023.

\bibitem{EH:2023}
B.~Ebner and N.~Henze.
\newblock On the eigenvalues associated with the limit null distribution of the
  {E}pps-{P}ulley test of normality.
\newblock {\em Statistical Papers}, 64(3):739--752, 2023.

\bibitem{ebner2021new}
B.~Ebner, S.~C. Liebenberg, and I.~J.~H. Visagie.
\newblock A new omnibus test of fit based on a characterization of the uniform
  distribution.
\newblock {\em Statistics}, 56(6):1364--1384, 2022.

\bibitem{FanEtAl2017}
Y.~Fan, P.~Lafaye~de Micheaux, S.~Penev, and D.~Salopek.
\newblock Multivariate nonparametric test of independence.
\newblock {\em Journal of Multivariate Analysis}, 153:189--210, 2017.

\bibitem{FM:2016}
G.~E. Fasshauer and M.~McCourt.
\newblock {\em Kernel-based approximation methods using MATLAB}.
\newblock Interdisciplinary mathematical sciences ; 19. World Scientific, New
  Jersey, 2016.

\bibitem{G:2004}
W.~Gautschi.
\newblock {\em Orthogonal polynomials : computation and approximation}.
\newblock Numerical mathematics and scientific computation. Oxford Univ. Press,
  Oxford, 2004.

\bibitem{HR:2020b}
E.~Hadjicosta and D.~Richards.
\newblock Integral transform methods in goodness-of-fit testing, {I}: the gamma
  distributions.
\newblock {\em Metrika}, 83(7):733--777, 2020.

\bibitem{HR:2020}
E.~Hadjicosta and D.~Richards.
\newblock Integral transform methods in goodness-of-fit testing, {II}: the
  {W}ishart distributions.
\newblock {\em Annals of the Institute of Statistical Mathematics},
  72(6):1317--1370, 2020.

\bibitem{H:1990}
N.~Henze.
\newblock An approximation to the limit distribution of the {E}pps-{P}ulley
  test statistic for normality.
\newblock {\em Metrika}, 37(1):7--18, 1990.

\bibitem{HJ:2019}
N.~Henze and M.~D. Jim\'{e}nez-Gamero.
\newblock A new class of tests for multinormality with i.i.d. and {GARCH} data
  based on the empirical moment generating function.
\newblock {\em TEST}, 28(2):499--521, 2019.

\bibitem{HJ:2021}
N.~Henze and M.~D. Jim\'{e}nez-Gamero.
\newblock A test for {G}aussianity in {H}ilbert spaces via the empirical
  characteristic functional.
\newblock {\em Scandinavian Journal of Statistics}, 48(2):406--428, 2021.

\bibitem{HJM:2019}
N.~Henze, M.~D. Jim\'{e}nez–Gamero, and S.~G. Meintanis.
\newblock Characterizations of multinormality and corresponding tests of fit,
  including for garch models.
\newblock {\em Econometric Theory}, 35(3):510–546, 2019.

\bibitem{HM:2019}
N.~Henze and C.~Mayer.
\newblock More good news on the {HKM} test for multivariate reflected symmetry
  about an unknown centre.
\newblock {\em Annals of the Institute of Statistical Mathematics},
  72:741–770, 2020.

\bibitem{HN:2000}
N.~Henze and Y.~Y. Nikitin.
\newblock A new approach to goodness-of-fit testing based on the integrated
  empirical process.
\newblock {\em Journal of Nonparametric Statistics}, 12(3):391--416, 2000.

\bibitem{HW:1997}
N.~Henze and T.~Wagner.
\newblock A new approach to the bhep tests for multivariate normality.
\newblock {\em Journal of Multivariate Analysis}, 62(1):1 -- 23, 1997.

\bibitem{HNNS:2019}
M.~Hu\v{s}kov\'{a}, N.~Neumeyer, T.~Niebuhr, and L.~Selk.
\newblock Specification testing in nonparametric {AR-ARCH} models.
\newblock {\em Scandinavian Journal of Statistics}, 46(1):26--58, 2019.

\bibitem{Imhof}
J.~P. Imhof.
\newblock Computing the distribution of quadratic forms in normal variables.
\newblock {\em Biometrika}, 48:419--426, 1961.

\bibitem{circular}
S.~R. Jammalamadaka, M.~D. Jim\'{e}nez-Gamero, and S.~G. Meintanis.
\newblock A class of goodness-of-fit tests for circular distributions based on
  trigonometric moments.
\newblock {\em SORT}, 43(2):317--336, 2019.

\bibitem{JAJB:2015}
M.~Jiménez-Gamero, M.~Alba-Fernández, P.~Jodrá, and I.~Barranco-Chamorro.
\newblock An approximation to the null distribution of a class of cramér–von
  mises statistics.
\newblock {\em Mathematics and Computers in Simulation}, 118:258--272, 2015.

\bibitem{K:1971}
R.~P. Kanwal.
\newblock {\em Linear integral equations : theory and technique}.
\newblock Academic Press, New York, 1971.

\bibitem{K:2001}
B.~Klar.
\newblock Goodness-of-fit tests for the exponential and the normal distribution
  based on the integrated distribution function.
\newblock {\em Annals of the Institute of Statistical Mathematics},
  53(2):338--353, 2001.

\bibitem{KG:2000}
V.~Koltchinskii and E.~Gin\'{e}.
\newblock Random matrix approximation of spectra of integral operators.
\newblock {\em Bernoulli}, 6(1):113--167, 2000.

\bibitem{Mardia&Jupp2000}
K.~Mardia and P.~Jupp.
\newblock {\em Directional statistics}.
\newblock Wiley Series in Probability and Statistics. John Wiley \& Sons, Ltd.,
  Chichester, 2000.

\bibitem{MT:2005}
M.~Matsui and A.~Takemura.
\newblock Empirical characteristic function approach to goodness-of-fit tests
  for the cauchy distribution with parameters estimated by mle or eise.
\newblock {\em Annals of the Institute of Statistical Mathematics},
  57(1):183--199, 2005.

\bibitem{meintanis2022bahadur}
S.~Meintanis, B.~Milo{\v{s}}evi{\'c}, and M.~Obradovi{\'c}.
\newblock Bahadur efficiency for certain goodness-of-fit tests based on the
  empirical characteristic function.
\newblock {\em Metrika}, pages 1--29, 2022.

\bibitem{M:2020}
B.~Milo\v{s}evi\'{c}.
\newblock Asymptotic efficiency of goodness-of-fit tests based on too–lin
  characterization.
\newblock {\em Communications in Statistics - Simulation and Computation},
  49(8):2082--2101, 2020.

\bibitem{milovsevic2021bahadur}
B.~Milo\v{s}evi\'c, {\relax Ya}.~{\relax Yu}. Nikitin, and M.~Obradovi\'c.
\newblock Bahadur efficiency of {EDF} based normality tests when parameters are
  estimated.
\newblock {\em Zapiski {N}auchnykh {S}eminarov POMI}, 501:203--217, 2021.

\bibitem{nikitinKnjiga}
{\relax Ya}.~{\relax Yu}. Nikitin.
\newblock {\em Asymptotic {E}fficiency of {N}onparametric {T}ests}.
\newblock Cambridge University Press, New York, 1995.

\bibitem{nikitinMetron}
{\relax Ya}.~{\relax Yu}. Nikitin and I.~Peaucelle.
\newblock Efficiency and local optimality of nonparametric tests based on {U}-
  and {V}-statistics.
\newblock {\em Metron}, 62(2):185--200, 2004.

\bibitem{NJ:2016}
F.~Novoa-Mu\~{n}oz and M.~D. Jim\'{e}nez-Gamero.
\newblock A goodness-of-fit test for the multivariate {P}oisson distribution.
\newblock {\em SORT}, 40(1):113--138, 2016.

\bibitem{NIST}
F.~W.~J. Olver, D.~W. Lozier, R.~F. Boisvert, and C.~W. Clark, editors.
\newblock {\em N{IST} handbook of mathematical functions}.
\newblock U.S. Department of Commerce, National Institute of Standards and
  Technology, Washington, DC; Cambridge University Press, Cambridge, 2010.

\bibitem{RW:2008}
C.~E. Rasmussen and C.~K.~I. Williams.
\newblock {\em Gaussian processes for machine learning}.
\newblock Adaptive computation and machine learning. MIT Press, Cambridge,
  Mass., 3. print edition, 2008.

\bibitem{RJ:2021}
G.-I. Rivas~Mart\'{\i}nez and M.-D. Jim\'{e}nez-Gamero.
\newblock Computationally efficient approximations for independence tests in
  non-parametric regression.
\newblock {\em Journal of Statistical Computation and Simulation},
  91(6):1134--1154, 2021.

\bibitem{RJ:2019}
G.~I. Rivas-Mart\'{\i}nez, M.~D. Jim\'{e}nez-Gamero, and J.~L. Moreno-Rebollo.
\newblock A two-sample test for the error distribution in nonparametric
  regression based on the characteristic function.
\newblock {\em Statistical Papers}, 60(4):1369--1395, 2019.

\bibitem{rublik1989optimality}
F.~Rubl{\'\i}k.
\newblock On optimality of the {LR} tests in the sense of exact slopes. {I}.
  general case.
\newblock {\em Kybernetika}, 25(1):13--14, 1989.

\bibitem{SW:1986}
G.~R. Shorack and J.~A. Wellner.
\newblock {\em Empirical Processes with Applications to Statistics}.
\newblock Wiley, 1986.

\bibitem{S:1976}
M.~A. Stephens.
\newblock {Asymptotic Results for Goodness-of-Fit Statistics with Unknown
  Parameters}.
\newblock {\em The Annals of Statistics}, 4(2):357 -- 369, 1976.

\bibitem{S:1977}
M.~A. Stephens.
\newblock Goodness of fit for the extreme value distribution.
\newblock {\em Biometrika}, 64(3):583--588, 1977.

\bibitem{WS:2003}
C.~Williams and J.~S. Shawe-Taylor.
\newblock The stability of kernel principal components analysis and its
  relation to the process eigenspectrum.
\newblock In S.~Becker, S.~Thrun, and K.~Obermayer, editors, {\em Advances in
  Neural Information Processing Systems 15}, pages 383--390. MIT Press, 2003.

\end{thebibliography}

\begin{appendix}
\section{Alternative eigenvalue approximation approaches}\label{app:AAA}

\subsection{Monte Carlo approach}\label{app:MCA}
This methodology was (not exclusively) introduced in \cite{EH:2023b} and \cite{FanEtAl2017}, and is associated with the quadrature method in the traditional numerical literature, see \cite{B:1977}, Chapter 3. It can also be found in machine learning theory, see \cite{RW:2008}, and for the approximation of the spectra of Hilbert-Schmidt operators, see \cite{KG:2000}. Assume that $\mu$ is a probability measure or equivalently that the weight function $w$ is scaled so that it is the corresponding probability density function. Then, let $Y\sim w$ be a random variable. Therefore, we can rewrite \eqref{int:eq} as
\begin{equation}\label{int:eq2}
\lambda f(s) = \E(K(s,Y)f(Y)),\quad s\in \R.
\end{equation}
An empirical counterpart to (\ref{int:eq2}) is found by letting $y_1,y_2,\ldots,y_N$, $N\in\N$, be independent realizations of $Y$ and approximating the expected value in (\ref{int:eq2}) by
\begin{equation}\label{int:eq3}
\E(K(s,Y)f(Y))\approx \frac1N\sum_{j=1}^N K(s,y_j)f(y_j),\quad s\in \R.
\end{equation}
Evaluate (\ref{int:eq3}) at the points $y_1,\ldots,y_n$ to get
\begin{equation}\label{int:eq3A}
\lambda f(y_i) = \frac1N\sum_{j=1}^N K(y_i,y_j)f(y_j),\quad i=1,\ldots,N,
\end{equation}
which is a system of $N$ linear equations. Writing $v=(f(y_1),\ldots,f(y_N))\in \R^N$ and $\widetilde{K}=(K(y_i,y_j)/N)_{i,j=1,\ldots,N}\in \R^{N\times N}$, we can rewrite (\ref{int:eq3A}) in matrix form
\begin{equation}\label{int:eq4}
\widetilde{K}v=\lambda v
\end{equation}
from which the (approximated) eigenvalues $\lambda_1,\ldots,\lambda_N$ can be computed explicitly. Note that for every eigenvalue $\lambda_j$ we have an eigenvector (say) $v_j\in\R^N$, whose components are the (approximated) values of the eigenfunctions (say) $f_j$ evaluated at $y_1,\ldots,y_N$. Obviously $\widetilde{K}$ is a random matrix, so that the calculated eigenvalues are random as well.

\subsection{Matrix-based operator' approximation approach}

In \cite{bozin2020new} the authors proposed the method for approximation of eigenvalues which is essentially based on the following steps: Recall \eqref{int:eq}, namely
\begin{equation*}
\lambda f(s) = \int_M K(s,t) f(t) \mu({\rm d}t),\quad s\in M,
\end{equation*}
with $\mu(\mbox{d}t)=w(t){\rm d}t$ and $\int_Mtw(t){\rm d}t<\infty$.

\begin{enumerate}
\item Consider the symmetrized operator

\begin{equation}\label{int:eq5}
\lambda f(s) = \int_M K(s,t) f(t) \sqrt{w(t)w(s)}{\rm d}t,\quad s\in M,
\end{equation}
which possesses the same spectrum as the one defined by \eqref{int:eq}.
\item  If $M=[-A,A]$, consider the sequence of symmetric linear operators defined by $(2m+1)\times(2m+1)$ matrices $M^{(m)}_{\omega}=||m^{(m)}_{i,j}||,\;|i|\leq m, |j|\leq m$
where
\begin{align*}
    m^{(m)}_{i,j}=K\left(\frac{iA}{m},\frac{jA}{m}\right)\sqrt{w\left(\frac{iA}{m}\right)w\left(\frac{jA}{m}\right)}.
\end{align*}
This sequence of operators converges in norm to the operator defined by \eqref{int:eq5}, and hence the spectra of these
two operators are at a distance that tends to zero.
\end{enumerate}
Observe that if $M=[a,b]$ for arbitrary $a<b$, the operator is suitably adjusted. In the case of unbounded support, step 2 is performed after truncating the operator \eqref{int:eq5} in such a manner that the resulting operator differs from the initial one on a set of negligible measure.

While this method is effective and straightforward to understand, one of its main disadvantages is that it may not demonstrate optimal computational efficiency, especially when dealing with complex kernel functions.

We demonstrate how this method operates using the kernels from Examples \ref{exm: exp} and \ref{exm: normality}.
The eigenvalues associated with the tests from \ref{exm: exp} and  \ref{exm: normality} are displayed in Tables \ref{tab:exm} and \ref{tab:norm} respectively.

\begin{table}[t]
\centering
\begin{tabular}{cclllll}
Kernel     & $\gamma$ &  $m=100$& $m=500$ & $m=1000$ & $m=2000$ & $m=5000$  \\ \hline
\eqref{kernel rho} & 0 &  2.70032e-1&2.71923e-1&2.72186e-1&2.7232e-1&2.72401e-1 \\
   & 1 & 1.00554e-1& 1.01133e-1 & 1.01227e-1& 1.01276e-1& 1.01306e-1\\
    &  2 &5.24446e-2& 5.26588e-2& 5.27049e-2&5.27296e-2& 5.27449e-2\\
     &  3& 3.22069e-2&3.22729e-2&3.22988e-2&3.23134e-2 &3.23227e-2 \\ \hline
     \eqref{K2.Laguerre} & 0 &1.93010e-1& 1.94551e-1&1.94746e-1& 1.94843e-1 &1.94901e-1 \\
   & 1 & 3.06451e-2& 3.08897e-2&3.09206e-2&3.09361e-2&3.09453e-2\\
    &   2 & 8.82674e-3& 8.89719e-3 &8.90608e-3& 8.91053e-3& 8.91320e-3\\
     &    3 & 3.43041e-3& 3.45777e-3& 3.46122e-3& 3.46295e-3& 3.46399e-3
\end{tabular}
\medskip

\caption{Approximation of the largest eigenvalue $\lambda_1(\gamma)$ from Example \ref{exm: exp} with $A=10$.} \label{tab:exm}
\end{table}

\begin{table}[t] \centering
\begin{tabular}{ccclllll}
Kernel     & $\gamma$ & $A$ & $m=100$& $m=500$ & $m=1000$ & $m=2000$ & $m=5000$  \\\hline
\eqref{K_Z} & 0.5 & 5 &1.19722 & 1.20201&1.20261&1.20291&1.20309\\
 &  1 & 4 & 4.49887e-1 & 4.51685e-1 & 4.51910e-1& 4.52023e-1 & 4.52091e-1\\
 &  2 & 3 & 1.40959e-1 & 1.41522e-1 & 1.41593e-1 & 1.41628e-1 & 1.41649e-1\\
 &  3 & 3 & 6.55722e-2 &6.58342e-2& 6.58671e-2&6.58836e-2&6.58934e-2\\ \hline
\eqref{ET} & 1.5 & 4 & 3.36331e-1 & 3.37675e-1 & 3.37843e-1 & 3.37928e-1 & 3.37979e-1\\
    &  2 & 3 & 8.72568e-2 & 8.76048e-2 & 8.76485e-2 &8.76704e-2& 8.76835e-2\\
& 3 & 3 & 1.38116e-2 & 1.38668e-2 & 1.38738e-2 & 1.38772e-2 & 1.38793e-2 \\ \hline
\eqref{K_Norm}& 0.5 & 5 &9.62013e-1 &9.65857e-1& 9.66340e-1 & 9.66582e-1 &9.66728e-1\\
    & 1 & 4 &5.97812e-1 &6.00201e-1  &6.00501e-1  &6.00651e-1  &6.00741e-1\\
     &  2 & 3 &3.96036e-1 & 3.97618e-1   & 3.97817e-1 &3.97916e-1&3.97976e-1 \\
& 3 & 3 &3.08925e-1   & 3.10160e-1& 3.10315e-1&3.10392e-1  &3.10439e-1
\end{tabular}
\medskip

\caption{Approximation of the largest eigenvalue $\lambda_1(\gamma)$ from Example \ref{exm: normality}} \label{tab:norm}
\end{table}

\section{Local Approximate Bahadur Efficiency}\label{ap: LABE}









In order to compute the approximate Bahadur efficiency, we need the following conditions to be met:
\begin{enumerate}
\item $T_n$ converges in distribution to a non-degenerate distribution function $F$,

\item $\log(1-F(t))=-\frac{a_Tt^2}{2}(1+o(1)),\; \to \infty$,

\item The limit in probability under an alternative from $\mathcal{G}$
$\lim_{n\rightarrow\infty}T_n/\sqrt{n}=b_T(\theta)>0
$ exists for $\theta \in \Theta_1$.
\end{enumerate}

Then, we have
\begin{equation}\label{BASlope}
c^{\ast}_T(\theta)=a_Tb_T^2(\theta)
\end{equation}

which represents the approximate Bahadur slope of $T_n$. Typically, under a close alternative where certain regularity conditions are met, i.e., when $\theta \to 0$, $b_T(\theta)$ can be expanded into a Maclaurin series. In the case of test statistics with an asymptotic normal distribution, the coefficient $a_T$ is equal to the inverse of the asymptotic variance. However, when the distribution is that of $\sum_{j=1}^{\infty}\lambda_jN_j^2,\;\lambda_1\geq\lambda_2\geq\cdot$, the coefficient $a_T$, which corresponds to $\sqrt{T_n}$, is equal to $\lambda_1^{-1}$ (for a detailed explanation, see \cite{cuparic2022}).

It's important to note that the local approximate Bahadur efficiency is an approximation of the local exact Bahadur efficiency, a measure that is significantly more complex to calculate in practice. While \cite{nikitinMetron} provided results for test statistics in the form of U- or V-statistics, the case involving unknown parameters still remains an open question.

\end{appendix}
\end{document}